\documentclass[reqno]{amsart}

\usepackage[top=1.1in, bottom=1in, left=1in, right=1in]{geometry}
\usepackage{amsfonts}
\usepackage{amssymb}
\usepackage{amsthm}
\usepackage{amsmath}
\usepackage{mathrsfs}
\usepackage{color}
\usepackage{enumerate}
\usepackage[numbers,sort&compress]{natbib}
\usepackage{pdfsync}
\usepackage{esint}
\usepackage{graphicx}
\usepackage{float}
\usepackage{caption}
\usepackage{subfigure}

\allowdisplaybreaks

\usepackage{tikz} 
\usetikzlibrary{calc}
\usetikzlibrary{intersections}
\usetikzlibrary{patterns}

\usepackage[colorlinks,
            linkcolor=black,
            anchorcolor=blue,
            citecolor=blue,
            ]{hyperref}

\makeatother

\numberwithin{equation}{section}

\newtheorem{theorem}{Theorem}[section]

\newtheorem{lemma}[theorem]{Lemma}
\newtheorem{remark}[theorem]{Remark}

\newtheorem{corollary}[theorem]{Corollary}
\numberwithin{equation}{section}

\def\({\left(}
\def\){\right)}

\def\calo{\mathcal{O}}

\def\bbe{{\mathbb{E}}}

\def\bbr{{\mathbb{R}}}
\def\bbp{{\mathbb{P}}}

\def\1{^{-1}}

\def\calo{{\mathcal{O}}}

\def\9{{\infty}}

\def\lbb{{\lambda}}
\def\a{{\alpha}}

\def\na{{\nabla}}

\def\wt{\widetilde}

\def\vf{{\varphi}}

\def\barr{\begin{array}}
\def\earr{\end{array}}

\def\({\left(}
\def\){\right)}
\def\<{\left<}
\def\>{\right>}

\def\wt{\widetilde}

\def\ol{\overline}
\def\ve{{\varepsilon}}

\begin{document}
	
\title[] {Recent progress on multi-bubble blow-ups
and multi-solitons to (stochastic)  focusing  nonlinear Schr\"odinger equations}

\author{Viorel Barbu}
\address{Octav Mayer Institute of
Mathematics (Romanian Academy)   and Al.I. Cuza University,
700506, Ia\c si, Romania.}
\email[Viorel Barbu]{vb41@uaic.ro}
\thanks{}

\author{Michael R\"ockner}
\address{Fakult\"at f\"ur
Mathematik, Universit\"at Bielefeld,  D-33501 Bielefeld, Germany.}
\email[Michael R\"ockner]{roeckner@math.uni-bielefeld.de}
\thanks{}

\author{Deng Zhang}
\address{School of Mathematical Sciences, CMA-Shanghai, Shanghai Jiao Tong University, China.}
\email[Deng Zhang]{dzhang@sjtu.edu.cn}
\thanks{}



\begin{abstract}
We review the recent progress on the long-time behavior for a
general class of focusing $L^2$-critical
nonlinear Schr\"odinger equations (NLS) with lower order perturbations.
Two canonical models are the stochastic NLS driven by
linear multiplicative noise and the classical deterministic NLS.
We show the construction and uniqueness of
the corresponding blow-up solutions
and solitons,
including the multi-bubble Bourgain-Wang type blow-up solutions
and non-pure multi-solitons,
which provide new examples for the mass quantization conjecture
and the soliton resolution conjecture.
The refined uniqueness of pure multi-bubble blow-ups and pure multi-solitons to NLS
under very low asymptotical rate is also reviewed.
Finally, as a new result,
we prove the qualitative properties of stochastic blow-up solutions,
including the concentration of mass,
universality of critical mass blow-up profiles,
as well as the vanishing of the virial at the blow-up time.
\end{abstract}

\maketitle

{
}

\section{Introduction}

We are concerned with a general class of focusing $L^2$-critical nonlinear Schr\"{o}dinger equations
with lower order perturbations
\begin{align} \label{gNLS}
    i\partial_t v +\Delta v + a_1 \cdot \nabla v + a_0 v +|v|^{\frac{4}{d}}v  =0 \tag{gNLS}
\end{align}
on $\bbr^d$, the coefficients of lower order perturbations are of form
\begin{align}
 & a_1(t,x)= 2 i \sum\limits_{l=1}^N \na \phi_l(x)h_l(t),   \label{a1-loworder} \\
 & a_0(t,x)= - \sum\limits_{j=1}^d \(\sum\limits_{l=1}^N \partial_j \phi_l(x)h_l(t)\)^2
          + i \sum\limits_{l=1}^N \Delta \phi_l(x) h_l(t), \label{a0-loworder}
\end{align}
with spatial functions $\phi_l \in C_b^\9(\bbr^d, \bbr)$,
and temporal functions $h_l \in C(\bbr^+; \bbr)$, $1\leq l\leq N<\infty$.

Equation \eqref{gNLS} is mainly motivated by the following two canonical models:

{\bf $\bullet$ Stochastic nonlinear Schr\"odinger equations.}
One important model closely related to \eqref{gNLS}
is the stochastic nonlinear Schr\"odinger equation (SNLS for short)
driven by linear multiplicative noise
\begin{align} \label{SNLS}
    & idX + \Delta X dt + |X|^{\frac 4d} X dt = - i \mu X dt +i X dW(t),  \tag{SNLS}
\end{align}
where
$W$ is a Wiener process
\begin{align}   \label{W-BM}
W(t,x)=\sum_{l=1}^N i\phi_l(x)B_l(t),\ \ x\in \bbr^d,\ t\geq 0,
\end{align}
$\{\phi_l\} \subseteq C_b^\9(\bbr^d, \bbr)$,
$\{B_l\}$ are independent standard $N$-dimensional real valued Brownian motions
on a normal stochastic basis $(\Omega, \mathscr{F}, \{\mathscr{F}_t\}, \bbp)$,
and
$\mu= \frac 12 \sum_{l=1}^N  \phi_l^2.$
The last term $XdW(t)$ in \eqref{SNLS} is taken in the sense of Gubinelli's controlled rough path,
which coincides with the usual It\^o stochastic integration
if the corresponding processes are $\{\mathscr{F}_t\}$-adapted
(cf. \cite{G04,FH14}).
We take $N<\infty$ for simplicity only.

The relationship between \eqref{gNLS} and \eqref{SNLS} can be seen
from the Doss-Sussman type transformation
\begin{align*}
    v:= e^{-W} X,
\end{align*}
which transforms the stochastic equation \eqref{SNLS} to
a random equation \eqref{gNLS}
with the temporal functions
$\{h_l\}$ being exactly the Brownian motions $\{B_l\}$.

The physical significance of SNLS is well known.
One significant model
arises from molecular aggregates
with thermal fluctuations,
and the multiplicative noise corresponds to
scattering of exciton by phonons,
due to thermal vibrations of the molecules.
Many physical and numerical experiments have been made to
study the noise effects on blow-up and solitons,
see, e.g., \cite{BCIR94,BCIRG95}
for noise effect on the coherence of the ground state solitary solution
in dimension two,
\cite{RGBC95} for the case of the critical quintic nonlinearity in dimension one.
We also   refer to \cite{BG09} for applications to open quantum systems,
where the martingale feature of mass plays an important role.

{\bf $\bullet$ Nonlinear Schr\"odinger equations.}
When $a_1,a_0=0$, i.e., the lower order perturbations vanish,
\eqref{gNLS} reduces to the canonical
focusing $L^2$-critical nonlinear Schr\"odinger equation (NLS for short)
\begin{align} \label{NLS}
    i\partial_t v +\Delta v +|v|^{\frac{4}{d}}v =0.   \tag{NLS}
\end{align}
NLS is one of fundamental dispersive PDEs and
is important in continuum mechanics, plasma physics and optics (\cite{DNPZ92}).
In particular, for the cubic nonlinearity in the critical dimension two,
the phenomenon of mass concentration near collapse gives a rigorous basis
to the physical concept of ``strong collapse'' (\cite{SS99}).
\eqref{NLS} has the conservation laws for the mass
\begin{align}   \label{mass-def}
   M(v):=  \|v\|_{L^2},
\end{align}
and the Hamiltonian
\begin{align}   \label{Hamil-def}
   H(v): = \frac 12 \|\na v\|_{L^2}^2  -  \frac{d}{2d+4} \|v\|_{L^{2+\frac 4d}}^{2+\frac 4d}.
\end{align}

The $H^1$ local well-posedness of \eqref{NLS} is well-known,
see, e.g., \cite{C03, T06}.
The SNLS has been studied in \cite{BD03,BM14,BD99},
by using stochastic Strichartz estimates.
Regarding the pathwise solvability of \eqref{SNLS} and the related \eqref{gNLS},
we refer to \cite{BRZ14,BRZ16,BRZ16.1,HRZ18,Z17,Z18} and the related
applications in control theory \cite{BRZ18,Z20}.
Let us mention that,
when treating the stochastic equation \eqref{SNLS},
the Doss-Sussman transformation (or, rescaling transformation)
enables us to perform sharp analysis in a pathwise way,
which is in general difficult for the usual It\^o calculus.
This pathwise treatment, particularly,
gives path-by-path uniqueness and hence the cocycle property,
i.e., we obtain a random dynamical system
(cf. \cite{CF94}).
This provides a convenient platform
to analyze the long-time dynamical mechanism
including the current blow-up and soliton dynamics.
For the interested readers,
we refer to \cite{Z18.2} for a review of the rescaling approach.

The long-time behavior of solutions to \eqref{gNLS} is much more delicate.
An important role here is played by the ground state,
which is the unique positive radial solution to the nonlinear elliptic equation
\begin{align} \label{equa-Q}
    \Delta Q - Q + Q^{1+\frac{4}{d}} =0.
\end{align}
In view of \cite{W83,D15},
the mass of the ground state is the sharp threshold for the
global well-posedness, scattering and blow-up.
We also would like to mention that,
for the case of energy-critical NLS,
Kenig and Merle \cite{KM06} developed the concentration-compactness and rigidity method,
and proved a sharp condition,
characterized by the Aubin-Talenti solution, for the global well-posedness, scattering
and blow-up.

In the critical mass case where $\|v\|_{L^2} = \|Q\|_{L^2}$,
the following two important dynamics arise:
\begin{enumerate}
  \item[$\bullet$] {Pseudo-conformal blow-up solutions}
\begin{align}  \label{S-blowup-intro}
    S_T(t,x)=(w(T-t))^{-\frac d2}Q \(\frac{x-x^*}{w(T-t)}\)
             e^{ - \frac i 4 \frac{|x-x^*|^2}{T-t} + \frac{i}{w^2(T-t)} + i\vartheta},
\end{align}
  \item[$\bullet$] {Solitary wave}
\begin{align}  \label{W-soliton-intro}
  W(t,x):=w^{-\frac d2}Q \(\frac{x- c t}{w} \)e^{i(\frac 12 c \cdot x-\frac{1}{4}|c|^2t+w^{-2}t+\vartheta)},
\end{align}
\end{enumerate}
where $x^*,c\in \bbr^d$,  $w>0$ and $T,\vartheta \in \bbr$.
Both dynamics are related closely to each other
through the
{\it pseudo-conformal transform}
in the pseudo-conformal space
$\Sigma :=\{u\in H^1: xu\in L^2\}$:
\begin{align} \label{pseu-conf-transf}
   S_T(t,x) = \mathcal{C}_T(W)(t,x):= \frac{1}{(T-t)^{\frac d2}} W \(\frac{1}{T-t}, \frac{x}{T-t}\) e^{-i\frac{|x|^2}{4(T-t)}}, \ \ t\not =T,\ x^* =c.
\end{align}
Note that $S_T$ blows up at time $T$,
and $x^*$ is the singularity
corresponding to the velocity $c$ of $W$.
In the seminal paper \cite{M93},
Merle proved that the pseudo-conformal blow-up solution is the unique critical mass blow-up solution
to the $L^2$-critical NLS, up to the symmetries of the equation.
It was conjectured that the only non-scattering solutions to \eqref{NLS},
up to symmetries,
are the soliton or the pseudo-conformal transformation of the soliton.
This conjecture has been recently solved by Dodson \cite{D21.1,D21.2}.

Moreover, in the small supercritical mass case where
$\|Q\|_{L^2} < \|v_0\|_{L^2} <\|Q\|_{L^2} + \ve$,
with small $\ve(>0)$,
one stable dynamics is the log-log blow-up solution,
which was first constructed by Perelman \cite{P01},
and has been extensively studied in a series of works
by Merle and Rapha\"el \cite{MR03,MR04,MR05.2,MR06}.

Another important dynamics is the Bourgain-Wang blow-up solution,
which was first constructed by Bourgain and Wang \cite{BW97},
and behaves asymptotically as a sum of a
singular profile $S_T$ and a regular profile $z$, i.e.,
\begin{align*}
   v(t) -S_T(t) -z(t) \to 0,\ \ as\ t\ \to T.
\end{align*}
Unlike the log-log blow-up solutions,
Bourgain-Wang solutions are unstable under $H^1$ perturbation
and lie on the boundary of two $H^1$ open sets of
global scattering solutions and log-log blow-up solutions (\cite{MRS13}).

In the even large mass regime,
the complete characterization of  general blow-up solutions
to $L^2$-critical NLS remains open.
In \cite{B00}, Bourgain raised an open problem on the
quantization property of blow-up solutions,
namely, whether the concentration of mass is of the form
$k\|Q\|_{L^2}^2$, $k\in \mathbb{Z}_+$.
See also \cite{BW97}.
Merle and Rapha\"el \cite{MR05} formulated precisely
the {\it mass quantization conjecture}:
blow-up solutions are expected to decompose into a singular part and
an $L^2$ residual,
and the singular part expands asymptotically as
multiple bubbles, each of which concentrates a mass of
no less than $\|Q\|_{L^2}^2$ at the blow-up point.
In particular,
the Bourgain-Wang solutions provide examples
of the mass quantization conjecture in the single-bubble case.
Multi-bubble blow-up solutions without regular profiles to \eqref{NLS} were first constructed
by Merle \cite{M90}.
Thus, a natural question is whether
there exist multi-bubble Bourgain-Wang
type blow-up solutions.

Furthermore,
according to the famous {\it soliton resolution conjecture},
global solutions to a nonlinear dispersive equation
are expected to decompose at large time
as a sum of solitons plus a scattering remainder.
Important progress has recently been made for the energy critical wave equation,
we refer to \cite{CKLS18,DJKM17,DKM13,DKM19}
and  references therein.
There are also many papers on the construction of multi-solitons
behaving like a sum of multi-solitons.
However,
the existence of {\it non-pure multi-solitons} (including dispersive part)
predicted by the soliton resolution conjecture
seemed not available,
see the lecture notes by Cazenave \cite{C20}.
An interesting question is then to construct this kind of {non-pure multi-solitons}.

Let us also mention that
the {\it uniqueness} of multi-solitons to $L^2$-(sub)critical NLS remains open
(see \cite{M18}).
The complete study of the {uniqueness problem} of multi-solitons
was done for the $L^2$-(sub)critical gKdV in the remarkable paper by Martel \cite{Ma05}.
Multi-solitons to $L^2$-supercritical gKdV were classified by C\^ombet \cite{Co11}.
We also would like to refer to \cite{JKL} for the classification of the strongly interacting kink-antikink pair.
For the NLS, recent progress has been made by C\^ote and Friederich \cite{CF20}
on the smoothness and  uniqueness in the (sub)critical case
when the asymptotic rate is large enough.
A natural question is thus to prove the uniqueness of solitons
particularly in the low asymptotic regime.

Turning to the stochastic case,
there are several physical and numerical papers on the study of
blow-up and solitons.
It was observed in \cite{BDM02,DL02,DL02.2} that
multiplicative noise
has the effect to delay blow-up,
while white noise even can prevent blow-up.
Moreover,
de Bouard and Debussche \cite{BD05} proved that, in the $L^2$-supercritical case,
the conservative noise  accelerates blow-up with positive probability.
Similar noise effect was also proved
for additive noise (\cite{BD02}).
In  \cite{BRZ14.3},
the authors proved that the non-conservative noise can provide damping effects,
and so prevent explosion with high probability.
Recently,
in \cite{MRRY20,MRY20},
quantitative behavior of blow-up solutions,
including the blow-up rate and blow-up profiles,
was studied by numerical experiments.
The main challenges in the stochastic case
are the absence of pseudo-conformal symmetry
and even of a conservation law for the Hamiltonian,
due to the presence of lower order perturbations (or noise).

In this note,
we review the recent progress \cite{FSZ20,SZ19,SZ20,RSZ21,RSZ21.2} on the construction and uniqueness
of blow-up solutions and of solitons
for the nonlinear Schr\"odinger equations with lower order perturbations \eqref{gNLS}.
In particular, the construction of multi-bubble
Bourgain-Wang type blow-up solutions
and non-pure multi-solitons provide new examples for the mass quantization conjecture
and the soliton resolution conjecture.
Furthermore, we also review the refined uniqueness results in \cite{CSZ21}
for pure multi-solitons to the $L^2$-critical \eqref{NLS}
in the low asymptotic rate regime.
Finally,
inspired by \cite{Ba,HK05,M93,M98},
as a new result,
in \S \ref{Sec-Mass-Concent-Univ} we prove qualitative properties of stochastic blow-up solutions,
including the concentration of mass, the universality of blow-up profiles with critical mass
and the vanishing of the virial at the blow-up time.

\section{Multi-bubble blow-up solutions and multi-solitons} \label{Sec-Multi-Blow-Soliton}

Let us first review the stochastic single-bubble blow-up solutions
in both the critical and small supercritical mass cases.
Then, we show the existence and uniqueness results
for the multi-bubble blow-up solutions and multi-solitons
in \S \ref{Subsec-Multi-BW} and \S \ref{Subsec-Soliton}, respectively.
In \S \ref{Subsec-Uniq},
we show the refined uniqueness results in the low asymptotical rate regime.

\subsection{Single bubble blow-up solutions}  \label{Subsec-single-bubble}
For the spatial functions $\{\phi_l\}$ we assume:
\begin{enumerate}
   \item[(A0)] Asymptotical flatness:
  For every $1\leq l\leq N$ and multi-index $\nu \not = 0$,
\begin{align*}
   \lim_{|x|\to \9} \<x\>^2 |\partial_x^\nu \phi_l(x)| =0.
\end{align*}

  \item[(A1)] Flatness at the origin:
  For every $1\leq l\leq N$ and multi-index $0\leq|\nu|\leq 5$,
\begin{align*}
   \partial_x^\nu \phi_l(0)=0.
\end{align*}
\end{enumerate}

The large spatial assumption in $(A0)$ ensures the local well-posedness of \eqref{gNLS}
(cf. \cite{BRZ14,BRZ16}),
while the local spatial assumption in $(A1)$
is mainly used for the blow-up analysis.

\begin{theorem} [Critical mass blow-up solution to SNLS, \cite{SZ19}]  \label{Thm-mini-gNLS}
Let $d=1,2$.
Assume Assumptions $(A0)$ and $(A1)$ to hold.
Then,
there exists $\tau^* \in (0,\9)$
such that for any $T\in (0,\tau^*]$,
there exists $v_0\in H^1$
satisfying
$\|v_0\|_{L^2} = \|Q\|_{L^2}$
and the corresponding  solution $v$ to \eqref{gNLS}
blows up at time $T$.
Moreover,  there exist $\delta, C(T)>0$ such that for $t$ close to $T$,
\begin{align*}
     \|v(t) - S_T(t) \|_{H^1} \leq C(T)(T-t)^\delta.
\end{align*}
\end{theorem}

\begin{remark}
$(i)$ A direct application of Theorem \ref{Thm-mini-gNLS} gives the existence
of critical mass blow-up solutions to \eqref{SNLS}.
It was known that
solutions to \eqref{SNLS}  exist globally in the $H^1$-subcritical case
or for initial data below the ground state, see \cite{BD03,BM14,BRZ16,FX18.1,M19}.
Hence, the mass of the ground state still serves as a threshold
for the global well-posedness and blow-up
in the stochastic case.

$(ii)$ The proof of Theorem \ref{Thm-mini-gNLS} is based on the modulation method
developed in \cite{R-S}, including the geometric decomposition,
a bootstrap device and the backward propagation from the singularity.
Delicate analysis has also been made to control the variation of the Hamiltonian,
related estimates of geometric parameters and the polynomial type of perturbation orders
when linearizing \eqref{gNLS} around the ground state.

$(iii)$
We also mention that the control of the first order term in \eqref{gNLS} is non-trivial
and is based on the local smoothing estimates (see, e.g., \cite{MMT08,Z17}).
The local smoothing estimates also play the key role in the well-posedness
of quasi-linear Schr\"odinger equations,
see \cite{KPV04,KPRV06} and the references therein.
\end{remark}

The next result is concerned with the construction of log-log blow-up solutions
to \eqref{gNLS}
in the small supercritical mass regime.
We refer to  \cite{P01,F17,MR03,MR04,MR05.2,MR06}
for \eqref{NLS} in the deterministic case.
For simplicity,
we focus on the case where $\{\phi_l\}$ are Schwartz functions.

\begin{theorem} [Log-log blow-up solution to SNLS, \cite{FSZ20}]  \label{Thm-loglog-SNLS}
Consider \eqref{gNLS}
with $d=1,2$ and $\{\phi_l\}$ being Schwartz functions.
Then, there exists an initial datum $v_0 \in H^1$
such that
the corresponding solution $v$ to \eqref{gNLS}
blows up in finite time $T$ according to a log-log law in the sense that
there exist parameters $(x, \gamma, \lambda) \in C^1((0,T); \bbr^d \times \bbr \times \bbr^+)$,
such that
\begin{equation*}
   v(t,x)=\frac{1}{\lambda^{d/2}(t)}(Q+\epsilon)(t,\frac{x-x(t)}{\lambda(t)})e^{i\gamma(t)} ,\ \ t\in (0,T),
\end{equation*}
with
\begin{equation*}
\lambda(t)^{-1}\sim \sqrt{\frac{\ln |\ln T-t|}{T-t}},
\text{ and } \int |\nabla \epsilon|^{2}+|\epsilon|^{2}e^{-|y|} dy \xrightarrow {t\rightarrow T} 0,
\end{equation*}
and the blow up point $x(t)$ converges as $t\to T$.
\end{theorem}

\subsection{Multi-bubble Bourgain-Wang type blow-up solutions}  \label{Subsec-Multi-BW}

This subsection is concerned with the multi-bubble blow-up solutions to \eqref{gNLS}
at the given distinct singularities $\{x_k\}$.
Similarly to $(A0)$ and $(A1)$,
we assume the following conditions for the spatial functions $\{\phi_l\}$
and the regular profiles $z^*$ (see, e.g., \cite{MR05,MRS13}).

\begin{enumerate}
   \item[$(H1_{\nu_*})$] {Asymptotical flatness}:
  For any multi-index $\nu \not = 0$ and $1\leq l\leq N$,
\begin{align*}
   \lim_{|x|\to \9} \<x\>^2 |\partial_x^\nu  \phi_l(x)| =0.
\end{align*}
  {Flatness at singularities}:
  There exists $\nu_*\in \mathbb{N}^+$ such that
  for every $1\leq l\leq N$ and every multi-index $\nu$ with  $|\nu|\leq \nu_*$,
\begin{align*}
   \partial_x^\nu  \phi_l(x_k)=0,\ \ 1\leq l\leq N,\ 1\leq k\leq K.
\end{align*}
  \item[$(H2_{m,\alpha^*})$] {Smallness}:
  There exist $m\in \mathbb{N}$, $\alpha^*\in (0,\infty)$
  such that for every regular profile $z^*$ satisfying
  \begin{align}  \label{z*-norm-small}
    & \|z^*\|_{H^{2m+2+d}}  \leq \a^*, \ \
      \|\<x\> z^*\|_{H^1}  \leq \a^*,
  \end{align}
  we have:

  {Flatness at singularities}:
  For any multi-index $\nu$ with $|\nu|\leq 2m$,
  \begin{align*}
    & \partial_x^\nu z^*(x_k) =0, \ \ 1\leq k\leq K.
  \end{align*}
\end{enumerate}

We first show the existence and conditional uniqueness of pure multi-bubble
blow-up solutions without the regular profile.

\begin{theorem} [Pure multi-bubble blow-up solutions to gNLS, \cite{SZ20}]\label{Thm-Multi-Blowup}
Consider \eqref{gNLS} with $d=1,2$.
Assume that $(H1_{\nu_*})$ holds with $\nu_*\geq 5$.
Let $K\in \mathbb{N}^+$, $T\in \bbr^+$, $\{\vartheta_j\}_{j=1}^K \subseteq \bbr$,
$\zeta\in (0,1)$.
Then:

(a) For any distinct points $\{x_k\}_{k=1}^K \subseteq \mathbb{R}^d$,
$w\in (0,\infty)$,
there exists $\ve>0$ such that for any $\{w_k\}_{k=1}^K \subseteq \mathbb{R}^+$
with $\max_{1\leq k\leq K}|w_k-w|\leq \ve$,
assertions $(i)$ and $(ii)$ below hold.

(b) For any  $\{w_k\}_{k=1}^K \subseteq \mathbb{R}^+$ there exists $\ve>0$
such that for any distinct points $\{x_k\}_{k=1}^K \subseteq \bbr^d$
with $\min_{1\leq k\leq K}|x_j-x_k|\geq \ve^{-1}$
assertions $(i)$ and $(ii)$ below hold,
where:

{$(i)$ Existence:}
There exists $\tau^*>0$ small enough
such that for any $T\in (0,\tau^*)$,
there exist $v_0\in \Sigma$ and a corresponding blow-up solution $v$ to \eqref{gNLS}
satisfying
\begin{align*}
\|v(t)-\sum_{k=1}^KS_k(t)\|_{\Sigma}\leq C(T-t)^{\frac{1}{2}(\nu_*-5)+ \wt \zeta},\ \ t\in [0,T),
\end{align*}
where $C>0$, $\wt \zeta\in (0, 1)$,
and $\{S_k\}$ are the pseudo-conformal blow-up solutions
\begin{align}   \label{Sj-blowup}
S_{k}(t,x)=(w_k(T-t))^{-\frac d2}Q \(\frac{x-x_k}{w_k(T-t)}\)
             e^{ - \frac i 4 \frac{|x-x_k|^2}{T-t} + \frac{i}{w_k^2(T-t)} + i\vartheta_k}.
\end{align}

{$(ii)$ Conditional uniqueness:}
There exists a unique blow-up solution $v$ to \eqref{gNLS}
satisfying
\begin{align*}
\| v(t) -\sum_{k=1}^KS_k(t)\|_{H^1}\leq C(T-t)^{3+\zeta},\ \ t\in[0,T),
\end{align*}
where $C>0$.
\end{theorem}

\begin{remark}
$(i)$ The above uniqueness result shows that multi-bubble blow-up solutions
are unique in the energy class with the asymptotic rate $(T-t)^{3+}$.
Later in Theorem \ref{Thm-Uniq-Blowup}, we can enlarge the unique energy class
with a much lower convergence rate $(T-t)^{0+}$
in the case of \eqref{NLS}.

$(ii)$ The proof of Theorem \ref{Thm-Multi-Blowup}
is based on a localization procedure
and, in particular, the delicate decoupling between the interactions of different
remainder profiles in the geometric decomposition.
In order to respect the multi-bubble structure,
a new generalized energy is constructed,
which permits to decouple different bubbles and maintains the key
monotonicity property.
\end{remark}

Furthermore, for the multi-bubble Bourgain-Wang type blow-up solutions we have

\begin{theorem} [Multi-bubble Bourgain-Wang solutions to gNLS, \cite{RSZ21}] \label{Thm-BW-gNLS}
Consider \eqref{gNLS} with $d=1,2$.
Let $K\in \mathbb{N}^+$, $T\in \bbr^+$, $\{\vartheta_k\}_{k=1}^K \subseteq \bbr$,
$\zeta\in (0,1)$.
Then, assertions (a) and (b) from Theorem \ref{Thm-Multi-Blowup} hold
with (i), (ii) there, replaced by:

$(i)$ {\it Existence:}
If $(H1_{\nu_*})$ holds for $\nu_*\geq 5$
and $(H2_{m,\alpha^*})$ for some $\alpha^*>0$,
$m\geq 3$ if $d=2$, and $m\geq 4$ if $d=1$,
then there exists $\ve^*\in (0,\alpha^*]$
such that for every regular profile $z^*$
satisfying \eqref{z*-norm-small}
with $\alpha^*$ replaced by $\ve^*$,
there exists a solution $v$ to \eqref{gNLS}
satisfying that for $t$ close to $T$,
\begin{align*}
&\|v(t) -\sum_{k=1}^KS_k(t) - z(t) \|_{L^2}
 + (T-t) \|v(t) -\sum_{k=1}^KS_k(t) - z(t) \|_{\Sigma} \leq C(T-t)^{\frac {1}{2}(\kappa -1)},
\end{align*}
where  $\kappa:= (m+\frac d2-1)\wedge (\upsilon_*-2)$,
$C>0$,
$\{S_k\}$ are as in \eqref{Sj-blowup},
and $z$ is the unique solution of the equation
\begin{align*}
 & i \partial_t z + \Delta z + a_1 \cdot \nabla z+a_0 z +|z|^{\frac 4d} z =0,   \\
    & z(T) = z^*,  \notag
\end{align*}
where the coefficients $a_1, a_0$ are given by \eqref{a1-loworder} and \eqref{a0-loworder}, respectively.

$(ii)$ {\it Conditional uniqueness:}
If $(H1_{\nu_*})$ holds with $\nu_* \geq 12$
and $(H2_{m,\alpha^*})$ holds for some $\alpha^*>0$
and $m\geq 10$,
then
there exists a unique solution $v$ to \eqref{gNLS} satisfying
that for $t$ close to $T$,
\begin{align*}
\|v(t) -\sum_{k=1}^KS_k(t) - z(t) \|_{L^2}
+ (T-t)\|\na v(t) -\na \sum_{k=1}^KS_k(t) - \na z(t) \|_{L^2}\leq C(T-t)^{4+\zeta}.
\end{align*}
\end{theorem}

\begin{remark}
$(i)$ To the best of our knowledge,
Theorem \ref{Thm-BW-gNLS} provides the first examples of multi-bubble blow-ups
with a regular profile.
As a direct application, it provides new examples
for the mass quantization conjecture,
the constructed solution satisfies that as $t\to T$,
   \begin{align*}
       & |v(t)|^2 \rightharpoonup \sum\limits_{k=1}^K \|Q\|_{L^2}^2 \delta_{x_k} + |z^*|^2, \\
       &  v(t) \to z^*\ \ in\ L^2\(\bbr^d - \bigcup\limits_{k=1}^K B(x_k, R)\)
   \end{align*}
for any $R>0$.
In particular, it concentrates the mass $\|Q\|_{L^2}^2$ at each singularity
and the remaining part converges to a regular profile $z^*$.

$(ii)$ The conditional uniqueness result reveals the rigidity of the flow around
multi-bubble pseudo-conformal blow-up solutions and the regular profile.

$(iii)$ The main effort of the proof is dedicated to decoupling
the interactions between three types of profiles:
the main blow-up profile, the regular profile and the remainder profile.
Unlike in \cite{MRS13},
because multi-bubble blow-up solutions are in general not radial,
two new geometric parameters are introduced,
which in turn require coercivity type control of the energy.
Let us also mention that
additional temporal regularity is gained
by subtracting the energy of the regular profile from that of the approximate solutions,
which is important to run the bootstrap arguments.
\end{remark}

\subsection{(Non-pure) multi-solitons}  \label{Subsec-Soliton}

By virtue of the pseudo-conformal invariance of $L^2$-critical NLS,
Theorem \ref{Thm-BW-gNLS} also gives the corresponding results
of non-pure multi-solitons to \eqref{NLS}.

\begin{theorem} [Non-pure multi-solitons to NLS, \cite{RSZ21}] \label{Thm-nonpure-soliton-gNLS}
Consider \eqref{NLS} with $d=1,2$.
Let $K\in \mathbb{N}^+$,  $\{\vartheta_k\}_{k=1}^K \subseteq \bbr$,
$\zeta\in (0,1)$.
Then, assertions (a) and (b) from Theorem \ref{Thm-Multi-Blowup} hold
with $x_k$ replaced by $c_k$, $1\leq k\leq K$,
and (i) and (ii) there, replaced by:

$(i)$ Existence:
If $(H2_{m,\alpha^*})$ holds for some $\alpha^*>0$ with $m\geq 6$,
then there exists $\ve^*\in (0,\alpha^*]$
such that for every regular profile $z^*$
satisfying \eqref{z*-norm-small}
with $\alpha^*$ replaced by $\ve^*$,
there exists a solution $u$ to \eqref{NLS} satisfying
\begin{align*}
   \|u(t) - \sum\limits_{k=1}^K W_k(t)-\widetilde{z}(t)\|_{\Sigma} \leq Ct^{-\frac 12 \kappa + \frac 52},\ \ for\ t\ large\ enough,
\end{align*}
where $\kappa = m+ \frac d2-1$, $C>0$,
$\{W_k\}$ are the solitary waves
\begin{align}  \label{Wj-soliton}
W_k(t,x) =w_k^{-\frac d2}Q \(\frac{x-c_k t}{w_k} \)e^{i(\frac 12 c_k\cdot x-\frac{1}{4}|c_k|^2t+ w_k^{-2}t+\vartheta_{k})},
\end{align}
and $\wt z$ corresponds to the regular part $z$ in Theorem \ref{Thm-BW-gNLS}
through the inverse of the pseudo-conformal transform:
\begin{align*}
   \wt z(t,x)  = \mathcal{C}_T^{-1} z (t,x )= t^{-\frac d2} z\(T-\frac 1t, \frac x t\) e^{i\frac{|x|^2}{4t}}.
\end{align*}

$(ii)$ Conditional uniqueness:
If $(H2_{m,\alpha^*})$ holds for some $\alpha^*>0$ with $m\geq 16$,
then
there exists a unique non-pure multi-soliton $u$ to \eqref{NLS}
satisfying
\begin{align*}
    \|u(t) - \sum\limits_{k=1}^K W_k(t)-\widetilde{z}(t)\|_{\Sigma} \leq Ct^{-5-\zeta},\ \ for\ t\ large\ enough.
\end{align*}

\end{theorem}

\begin{remark}
$(i)$ To our best knowledge,
Theorem \ref{Thm-nonpure-soliton-gNLS} provides first examples of non-pure multi-solitons
including a dispersive part to $L^2$-critical \eqref{NLS}
predicted by the soliton resolution conjecture.

$(ii)$ It is worth noting that
the uniqueness of non-pure multi-solitons holds in the energy class
with the decay rate $t^{-5-}$,
which is larger than the usual class of fast exponential decay rate,
which multi-solitons naturally lie in (see, e.g., \cite{LeLP15,LeT14}).
For multi-solitons without scattering part,
we obtain the uniqueness with even lower asymptotic rate,
see Theorems \ref{Thm-Uniq-Solitons-H1} and \ref{Thm-Uniq-Solitons-Sigma} below.
\end{remark}

Unlike the $L^2$-critical NLS,
because of the absence of pseudo-conformal invariance,
stochastic multi-solitons to \eqref{SNLS} cannot be obtained from
multi-bubble blow-ups in \S \ref{Subsec-Multi-BW}.
Hence, the construction of stochastic multi-solitons requires
a direct method at the level of solitons.
This is the content of Theorem \ref{Thm-Soliton-SNLS} below.

Below we consider the $L^2$-(sub)critical SNLS of similar structure as that of \eqref{SNLS}:
\begin{equation}  \label{SNLS-Sub}
  dX(t) = i\Delta X(t)dt + i|X(t)|^{p-1}X(t) dt - \mu(t) X(t) dt + \sum\limits_{k=1}^N X(t) G_k(t)dB_k(t),
\end{equation}
where $1<p\leq 1+\frac 4d$, $d\geq 1$, $T_0\geq 0$,
$\{B_k\}$ and $\mu$ as in \eqref{SNLS},
$G_k(t,x)=  i\phi_k(x)g_k(t),\ \ x\in \bbr^d,\ t\geq 0$,
$\{\phi_k\} \subseteq C_b^\9(\bbr^d, \bbr)$,
$\{g_k\}$ are $ \{\mathscr{F}_t\}$-adapted processes
with paths in $C^\a(\bbr^+, \bbr)$,
$\a\in(\frac 13, \frac 12)$,
that are controlled by $\{B_k\}$
and $X(t) G_k (t) dB_k(t)$ is taken in the sense of controlled rough paths
(see, e.g., \cite{G04,FH14}).

The basic conditions here are as follows:

{(H1)'}
  For every $1\leq l\leq N$,
  \begin{align*}
       \lim\limits_{|x|\to \9} |x|^2 |\partial_x^\nu \phi_l(x)|=0, \ \ \nu \not =0.
  \end{align*}

{(H2)'}
  For every $1\leq l\leq N$, $\{g_l\}$
  are $\{\mathscr{F}_t\}$-adapted continuous processes
  and controlled by the Brownian motions $\{B_l\}$,
  i.e., $\{g_l\} \subseteq \mathscr{D}_B^{2\a}(\bbr^+; \bbr^N)$
  with Gubinelli  derivative $\{g'_{lj}\}_{j,l=1}^N$.
  In addition, for $1\leq l\leq N$,
  $\phi_l$ and $g_l$ satisfy one of the following two cases:

  {\rm Case (I)':}
  $g_l\in L^2(\bbr^+)$, $\bbp-a.s$., and
  there exists $c_l>0$ such that
  \begin{align} \label{phil-exp-decay}
     \sum\limits_{|\nu|\leq 4} |\partial^\nu \phi_l(x)| \leq C e^{-c_l |x|}.
  \end{align}

  {\rm Case (II)':}
  $\bbp$-a.s., $g_l\in L^2(\bbr^+)$
  and there exists $c^*>0$ such that for $t$ large enough,
  \begin{align}   \label{gl-t2-decay}
     \int_t^\9 g_l^2 ds \log \(\int_t^\9 g_l^2 ds \)^{-1} \leq \frac {c^*}{t^2}.
  \end{align}
  \qquad \qquad \ \ \ \
  In addition, let $\nu_*\in \mathbb{N}$ be such that $\phi_l$ satisfies
  \begin{align} \label{phil-poly-decay}
     \sum\limits_{|\nu|\leq 4} |\partial^\nu \phi_l(x)| \leq C |x|^{-\nu_*}.
  \end{align}

We see that
the asymptotic behavior \eqref{gl-t2-decay} is closely related to the Levy H\"older continuity
of Brownian motion.
For simplicity,
we focus on the case where $c_l=1$, $1\leq l\leq N$,
and set
\begin{align*}
   \phi(x) := \left\{
               \begin{array}{ll}
                 e^{-|x|}, & \hbox{in {\rm Case (I)'};} \\
                 |x|^{-\nu_*}, & \hbox{in {\rm Case (II)'}.}
               \end{array}
             \right.
\end{align*}

\begin{theorem} [Multi-solitons to SNLS, \cite{RSZ21.2}] \label{Thm-Soliton-SNLS}
Consider \eqref{SNLS-Sub} with $1 <p \leq  1+\frac 4d$, $d\geq 1$.
Let $w_k^0>0$, $\theta_k^0 \in \bbr$,
$x_k^0 \in \bbr^d$,
$v_k \in \bbr^d\setminus\{0\}$, $1\leq k\leq K$, such that
$v_j\not = v_k$ for any $j\not = k$.
Assume $(H1)'$ and $(H2)'$ with $\nu_*$ sufficiently large in {\rm Case (II)'}.

Then, for $\bbp$-a.e. $\omega\in \Omega$,
there exist $T_0= T_0(\omega)$ sufficiently large
and $X_*(\omega)\in H^1$,
such that there exists an $H^1$ solution $X(\omega)$ to \eqref{SNLS-Sub} on $[T_0,\9)$
satisfying
$X(\omega, T_0) = X_*(\omega)$ and
\begin{align}  \label{X-Rk-asym}
     \|e^{-W_*(t)}X(t) - \sum\limits_{k=1}^K R_k(t) \|_{H^1} \leq C \int_{t}^{\infty}s\phi^{\frac 12}(\delta s)ds, \ \ t\geq T_0.
\end{align}
Here, $\{R_k\}$ are the solitary waves
\begin{align*}
    R_k(t,x):=Q_{w_{k}^0} \(x-v_k t-x_{k}^0 \)e^{i(\frac 12 v_k\cdot x-\frac{1}{4}|v_k|^2t+(w_{k}^0)^{-2}t+\theta_{k}^0)},
\end{align*}
with $Q_{w}(x) :=w^{-\frac {2}{p-1}}Q \(\frac{x}{w} \)$, and
\begin{align*}
  W_*(t,x) =- \sum\limits_{l=1}^N \int_t^\9 i \phi_l(x) g_l(s) d B_l(s).
\end{align*}

Moreover, in the $L^2$-subcritical case $1<p<1+\frac 4d$,
there exists a solution $X$ to \eqref{SNLS-Sub} on
the whole time interval $[0,\9)$,
satisfying the asymptotic behavior \eqref{X-Rk-asym}.
\end{theorem}

\begin{remark}
$(i)$ We note that
the temporal asymptotic rate in \eqref{X-Rk-asym}
can be of either exponential and polynomial type, respectively,
in {\rm Case (I)'} and  {\rm Case (II)'},
which is closely related to the spatial decay rate of the noise.
This reflects the effect of noise on the soliton dynamics.

$(ii)$ The proof of Theorem \ref{Thm-Soliton-SNLS}
is based on the rescaling approach,
involving two types of Doss-Sussman transforms
as in the stochastic scattering context \cite{HRZ18},
and on the modulation method in \cite{CF20,MMT06,MM06}.
One crucial ingredient is the monotonicity of the
Lyapunov type functional constructed in \cite{MMT06}.
Let us also mention that
the geometric decomposition in the proof is chosen in a quite unified manner
in both the $L^2$-subcritical and critical cases.
\end{remark}

\subsection{Refined uniqueness}  \label{Subsec-Uniq}

In this subsection, we review the refined uniqueness results on the multi-bubble blow-up solutions
and multi-solitons to \eqref{NLS},
particularly,
in the low asymptotic regime.

\begin{theorem} [Refined uniqueness of multi-bubble blow-ups to NLS, \cite{CSZ21}]   \label{Thm-Uniq-Blowup}
Consider \eqref{NLS} in dimensions $d=1,2$.
Let $T\in \bbr$, $K \in \mathbb{N}\setminus\{0\}$,
$\{\vartheta_k\} \subseteq \bbr$.
Then, for any $\zeta\in (0,1)$,
for any distinct points $\{x_k\}_{k=1}^K \subseteq \mathbb{R}^d$,
$w>0$ (resp. $\{w_k\}_{k=1}^K \subseteq \mathbb{R}^+$),
there exists $\ve>0$ such that for any $\{w_k\}_{k=1}^K \subseteq \mathbb{R}^+$
with $\max_{1\leq k\leq K}|w_k-w|\leq \ve$
(resp. $\{x_k\}_{k=1}^K \subseteq \bbr^d$
with $\min_{1\leq k\leq K}|x_j-x_k|\geq \ve^{-1}$),
there exists a unique multi-bubble blow-up solution $v$ to \eqref{NLS}
satisfying
\begin{align*}
    \| v(t) - \sum\limits_{k=1}^K S_k(t) \|_{L^2}
    + (T-t) \| \na v(t) - \sum\limits_{k=1}^K \na S_k(t) \|_{L^2} = o(1),\ \ as\ t\ close\ to\ T,
\end{align*}
and additionally
\begin{align*}
     \frac{1}{T-t} \int_t^T \frac{1}{T-s} \int_s^T \|v (r) - \sum\limits_{k=1}^K S_k(r)  \|_{H^1}^2 dr ds
     = \calo((T-t)^{\zeta}),
\end{align*}
where $S_k$, $1\leq k\leq K$, are the
pseudo-conformal blow-up solutions
given by \eqref{Sj-blowup}.
Moreover, the asymptotic rate of the unique solution $v$
can be enhanced to the exponentially fast rate in
the pseudo-conformal space. Namely,
there exists $\delta>0$ such that
\begin{align*}
   \|v(t) - \sum\limits_{k=1}^K S_k(t) \|_{\Sigma} = \calo(e^{-\frac{\delta}{T-t}}), \ \ for\ t\ close\ to\ T.
\end{align*}

In particular,
the above results hold for
the multi-bubble blow-up solutions $v$ to \eqref{NLS} such that
\begin{align*}
   \| v(t) - \sum\limits_{k=1}^K S_k(t) \|_{H^1}  = \calo( (T-t)^{\zeta}),\ \ for\ t\ close\ to\ T.
\end{align*}
\end{theorem}

Regarding  the pure multi-solitons to \eqref{NLS},
we have the following uniqueness result in the energy class.

\begin{theorem} [Refined uniqueness in $H^1$ of multi-solitons to NLS, \cite{CSZ21}]\label{Thm-Uniq-Solitons-H1}
Consider the situation as in Theorem \ref{Thm-Uniq-Blowup}.
Then, for any $\zeta\in (0,1)$,
for any distinct speeds $\{c_k\}_{k=1}^K \subseteq \mathbb{R}^d$,
$w>0$ (resp. $\{w_k\}_{k=1}^K \subseteq \mathbb{R}^+$),
there exists $\ve>0$ such that for any $\{w_k\}_{k=1}^K \subseteq \mathbb{R}^+$
with $\max_{1\leq k\leq K}|w_k-w|\leq \ve$
(resp. $\{c_k\}_{k=1}^K \subseteq \bbr^d$
with $\min_{1\leq k\leq K}|c_j-c_k|\geq \ve^{-1}$),
there exists a unique multi-soliton $u$ to \eqref{NLS}
satisfying
\begin{align*}
\|u(t)- \sum_{k=1}^K W_k(t)\|_{H^1} = \calo\(\frac{1}{t^{2+\zeta}}\), \ \ for\ t\ large\ enough,
\end{align*}
where $\{W_k\}$ are the solitons given by \eqref{Wj-soliton}.
Moreover,
the unique multi-soliton $u$ converges exponentially fast
in the  pseudo-conformal space,
i.e., for some $\delta>0$,
\begin{align*}
\|u(t)- \sum_{k=1}^K W_k(t)\|_{\Sigma} = \calo(e^{- \delta t}),  \ \ for\ t\ large\ enough.
\end{align*}
\end{theorem}

In the case of the pseudo-conformal space,
the uniqueness class can be further enlarged
in the even lower asymptotic regime.

\begin{theorem} [Refined uniqueness in $\Sigma$ of multi-solitons to NLS, \cite{CSZ21}]  \label{Thm-Uniq-Solitons-Sigma}
Consider the situation as in Theorem \ref{Thm-Uniq-Solitons-H1}.
Then, for any $\zeta\in (0,1)$,
for any distinct non-zero speeds $\{c_k\}_{k=1}^K \subseteq \mathbb{R}^d \setminus \{0\}$,
$w>0$ (resp. $\{w_k\}_{k=1}^K \subseteq \mathbb{R}^+$),
there exists $\ve>0$ such that for any $\{w_k\}_{k=1}^K \subseteq \mathbb{R}^+$
with $\max_{1\leq k\leq K}|w_k-w|\leq \ve$
(resp. $\{c_k\}_{k=1}^K \subseteq \bbr^d \setminus \{0\}$
with $\min_{1\leq k\leq K}|c_j-c_k|\geq \ve^{-1}$),
there exists a unique multi-soliton $u$ to \eqref{NLS}
satisfying
\begin{align*}
\|u(t)- \sum_{k=1}^K W_k(t)\|_{\Sigma} = \calo\(\frac{1}{t^{\frac 12+\zeta}}\), \ \ for\ t\ large\ enough.
\end{align*}
Moreover, the unique multi-soliton $u$ converges exponentially fast to $\sum_{k=1}^K W_k$ in $\Sigma$.
\end{theorem}

\begin{remark}
$(i)$ The uniqueness of multi-solitons to \eqref{NLS}
was first obtained by C\^ote and Friederich \cite{CF20}
in the $L^2$-(sub)critical and critical cases,
provided that the convergence rate is $(1/t)^{N}$ for $N$ large enough.
Theorem \ref{Thm-Uniq-Solitons-H1} shows that
the uniqueness class of multi-solitons to $L^2$-critical \eqref{NLS}
can be enlarged in the low convergence regime with rate $(1/t)^{2+}$.

$(ii)$
The proof proceeds in several upgradation steps,
each step requires the monotonicity of different functionals
appropriately constructed with respect to the multi-bubble structure.
More delicately,
the analysis of the functionals relies on suitable estimates
of the remainder and geometric parameters in the previous step.

Another main difficulty is the presence of the localized mass,
which is absent in the single-bubble case.
The key idea is to upgrade the localized mass and the remainder together
in the upgradation procedure.

The a priori low asymptotic rate also gives rise to a challenging problem to
identify the energy.
The important fact here is that the remainder exhibits dispersion in the energy space
along a sequence,
which leads to the energy quantization phenomenon
that is the key towards the derivation of the refined energy estimate.
\end{remark}

\section{Mass concentration and universality}   \label{Sec-Mass-Concent-Univ}

This section is devoted to the qualitative descriptions of
the dynamical properties of blow-up solutions to \eqref{SNLS},
including the concentration of $L^2$-norm,
the universality of critical mass blow-up solutions,
as well as the vanishing of the virial at the blow-up time.

We consider the $H^1$ blow-up solution $X$ to \eqref{SNLS} on $[0,\tau^*)$,
where $d\geq 1$, the stochastic integration is taken in the sense of It\^o
and $\tau^*(>0)$ denotes the maximal existing time.
We note that by the construction in \cite{BRZ16},
$\tau^*$ is an $\{\mathscr{F}_t\}$-stopping time
and $\tau^*>0$, $\mathbb{P}$-a.s..
Hence, below we focus on the case $\mathbb{P}(\tau^*<\infty)>0$.

The first result is concerned with the mass concentration property
along a time sequence.

\begin{theorem} [Mass concentration] \label{Thm-L2-Con}
Assume the asymptotic flat condition $(A0)$ to hold.
Assume additionally that
$\mathbb{P}(\tau^* <\9)>0$.
Then,
there exists a measurable set $\Omega_0 \subseteq \Omega$ with
full probability
such that for every $\omega \in \{\tau^*<\9\} \cap \Omega_0$,
there exist
$t_n(\omega) \in (0,\tau^*)$
and $y_n(\omega) \in \mathbb{R}^d$, $n \in \mathbb{N}$,
such that $t_n(\omega)\to \tau^*(\omega)$ and  for any $R>0$,
\begin{align} \label{esti-L2-conc}
   \liminf\limits_{n\to\9}
   \int\limits_{|x-y_n(\omega)|\leq R} |X(t_n(\omega),x)(\omega)|^2 dx
   \geq \|Q\|_{L^2}^2.
 \end{align}
\end{theorem}

\begin{remark}
For the deterministic \eqref{NLS},
the concentration of the $L^2$-norm was proved in \cite{HK05,MT90,T90}
along all sequence of times converging to the blow-up time.
The reason that Theorem \ref{Thm-L2-Con} holds for a time sequence
is mainly due to the failure of the conservation of the Hamiltonian in the stochastic case.
The It\^o evolution of the Hamiltonian actually
contains an extra stochastic integration
(see \eqref{H-evo} below),
which fluctuates in time and in general is not uniformly bounded.

The key point here
is to make use of the specific structure of the stochastic integration,
which can be controlled by using the Banica type inequality
together with the Burkholder-Davis-Gundy inequality.
Moreover,
the quantitative large time behavior,
i.e., the law of the iterated logarithm of a martingale, is also used
in the control of the asymptotic behavior of the Hamiltonian of rescaled solutions.

\end{remark}

In the critical mass case where $\|X_0\|_{L^2} = \|Q\|_{L^2}$,
we have the universality of stochastic critical mass blow-up profiles,
which extends the result in \cite{W86} to the stochastic case.

\begin{theorem} [Universality of critical mass blow-up profiles]  \label{Thm-Self}
Assume the situation as in Theorem \ref{Thm-L2-Con} to hold.
Assume additionally that
$\|X(0)\|_{L^2} = \|Q\|_{L^2}$.
Then,
there exists a measurable set $\Omega_0 \subseteq \Omega$ with
full probability
such that
for any $\omega \in \{\tau^*<\9\} \cap \Omega_0$,
there exist $t_n(\omega) \in (0,\tau^*)$,
$y_n(\omega) \in \mathbb{R}^d$ and $\theta(t_n(\omega))\in \mathbb{R}$, $n \in\mathbb{N}$, such that
\begin{equation} \label{Strong-X-Q}
     (\lbb_n(\omega))^{\frac d2} X(t_n(\omega), \lbb_n\cdot + y_n(\omega)) e^{i\theta(t_n(\omega))} \to Q\ \ in\ H^1,\ as\ n\to \9,
\end{equation}
where $\lbb_n(\omega):= \frac{\|\na Q\|_{L^2}}{\|\na X(t_n(\omega))\|_{L^2}}$, $n \geq 1$.
\end{theorem}

A direct application of Theorem \ref{Thm-Self}
yields that the critical mass blow-up solutions concentrate
the whole mass at the singularity.

\begin{corollary} \label{Cor-L2-conv}
Assume the conditions of Theorem \ref{Thm-Self} to hold
and let $\Omega_0, y_n$ and $t_n$ be as in Theorem \ref{Thm-Self}.
Then,
we have for every $\omega\in \{\tau^*<\infty\} \cap \Omega_0$,
\begin{align} \label{X2-delta}
   |X(t_n(\omega),x)|^2 dx - \|Q\|_{L^2}^2 \delta_{y_n(\omega)} \to
   0\ \ in\ distribution,\ as\ n\to \9.
\end{align}
\end{corollary}

\begin{remark} \label{Rem-con}
$(i)$ In \cite{M93},
Merle proved the rigidity of the critical mass blow-up solutions to \eqref{NLS},
which are unique up to the symmetries of the equation.
See also \cite{HK05,R13} for simplified proofs.
In the stochastic case,
this strong rigidity is still unclear.
The obstructions are due to the lack of pseudo-conformal symmetry,
and the failure of the energy conservation and the virial identity.
In the energy class with asymptotic rate $(T-t)^{3+}$,
in view of Theorem \ref{Thm-Multi-Blowup}
the uniqueness of stochastic critical mass blow-up solutions holds.

$(ii)$ A sufficient condition for the
finite time blow-up, i.e., $\bbp(\tau^*<\9)>0$,
is that
$H(X_0) < 0$ and
$\sum_{k=1}^N \|\na \phi_k\|_{L^{\9}}\leq \ve$
with $\ve$ small enough,
see \cite[Proposition 3.1]{BRZ14.3}.
In particular,
if $H(X_0)<0$
and  $\{\phi_l\}$ are constants
(i.e., the noise $W$ is spatially independent),
then  \eqref{SNLS} reduces to the classical NLS
\begin{align*}
   idy =  \Delta y dt + |y|^{\frac 4d}y dt
\end{align*}
with $y(0)=X_0$.
Thus, the explosion time $\tau^*$ is exactly the same as in the deterministic case.
\end{remark}

In the general case,
the numerical results in  \cite{BDM02,DL02,DL02.2} suggest that,
though the smooth multiplicative noise cannot prevent blow-up, i.e., $\tau^*<\infty$, $\bbp$-a.s.,
it has the effect to delay blow-up.
The rigorous proof remains unclear.
In the following we show that along certain sequence $t_n\to \tau^*$,
the virial will tend to zero,
which might be of use to understand the blow-up time.

\begin{theorem} [Vanishing of the virial at blow-up time]    \label{Thm-Var-0}
Consider \eqref{SNLS} with $X_0\in H^1$.
Assume $(A0)$ to hold.
Assume additionally that $\|X_0\|_{L^2}=\|Q\|_{L^2}$ and
$\mathbb{P}(\tau^*<\9)=1$.
Let  $t_n$ and $y_n$, $n\in \mathbb{N}$, be as in Theorem \ref{Thm-Self}.
Then, there exists a measurable set $\Omega_1$ with full probability
such that for any $\omega\in \Omega_1$,
there exists $y^*(\omega)\in\mathbb{R}^d$
such that
$y_n(\omega)\to y^*(\omega)$ and
\begin{align} \label{conv-x-x*}
     \lim\limits_{n\to \9} \int |x-y^*(\omega)|^2 |X(t_n,x)(\omega)|^2 dx =0.
\end{align}
Moreover,
we have the following lower bound of the blow-up rate along the sequence $\{t_n(\omega)\}$,
\begin{align} \label{blowup-rate}
    \|\na X(t_n)(\omega)\|_{L^2}^2 \geq C(\tau^*(\omega)) (\tau^*(\omega) - t_n(\omega))^{-2}.
\end{align}
\end{theorem}

The proofs of Theorems \ref{Thm-L2-Con}, \ref{Thm-Self} and \ref{Thm-Var-0}
are inspired by \cite{Ba,HK05,M93,R13}
and are contained in the subsequent
\S \ref{Subsec-Concentrate}, \S \ref{Subsec-Self} and \S \ref{Subsec-Explosion},
respectively.

\subsection{Concentration of $L^2$-norm} \label{Subsec-Concentrate}

Let us start with the compactness result in \cite{HK05}.

\begin{theorem} [Theorem 1.1, \cite{HK05}] \label{Thm-Compact}
Let $\{v_n\}_{n=1}^\9$ be a bounded family of $H^1$ such that
\begin{align}
   \limsup\limits_{n\to \9} \|\na v_n\|_{L^2} \leq M,  \ \
   \limsup\limits_{n\to \9} \|v_n\|_{L^{2+\frac 4d}} \geq m.
\end{align}
Then,
there exists $\{x_n\}_{n=1}^\9 \subseteq \bbr^d$ such that,
up to a subsequence,
\begin{align*}
    v_n (\cdot + x_n) \rightharpoonup V\ weakly,
\end{align*}
with
\begin{align}
  \|V\|_{L^2} \geq (\frac{d}{d+2})^{\frac d4} ( \frac{m^{\frac d2 +1}}{M^{\frac d2}}) \|Q\|_{L^2}.
\end{align}
\end{theorem}

\begin{lemma} [Sharp Gagliardo-Nirenberg inequality, \cite{W83}] \label{Lem-GN}
For any $v\in H^1$, we have
\begin{align} \label{GN}
     \|v\|^{2+\frac 4d}_{L^{2+\frac 4d}} \leq (\frac{d+2}{d}) (\frac{\|v\|_{L^2}}{\|Q\|_{L^2}})^{\frac 4d} \|\na v\|_{L^2}^2.
\end{align}
In particular, for the Hamiltonian  given by \eqref{Hamil-def},
we have
\begin{align} \label{H-GN}
   H(v) \geq \frac 12 \(1-(\frac{\|v\|_{L^2}}{\|Q\|_{L^2}})^{\frac 4d} \) \|\na v\|^2_{L^2},
\end{align}
and $H(v)\geq 0$ if $\|v\|_{L^2} \leq \|Q\|_{L^2}$.
\end{lemma}

The following result extends Banica's inequality in \cite{Ba} to the stochastic case.

\begin{lemma} [Banica type estimate]    \label{Lem-integ}
Let $X$ be the solution to \eqref{SNLS} on $[0,\tau^*)$.
Suppose that $\|X_0\|_{L^2}\leq \|Q\|_{L^2}$.
Then, $\bbp$-a.s. for any function $\phi \in W^{1,\9}$ and for any $t\in [0,\tau^*)$
we have
\begin{align} \label{esti-integ}
    \bigg| {\rm Im} \int X(t) \na \overline{X(t)} \cdot \na \phi dx \bigg|
     \leq  \(2H(X(t)) \int |X(t)\na \phi|^2 dx \)^{\frac 12}.
\end{align}
\end{lemma}

{\it \bf Proof.}
Thanks to the conservation law of mass,
we adapt the arguments as in the proof of \cite[Lemma 2.1]{Ba}.
On one hand, by
\eqref{H-GN} and the conservation of mass, for any $\a\in \bbr$,
\begin{align*}
   H(e^{i\a \phi} X(t))
    \geq& \frac 12 \(1-(\frac{\|e^{i\a \phi}X(t)\|_{L^2}}{\|Q\|_{L^2}})^{\frac{4}{d}}\) \|\na (e^{i\a \phi}X(t))\|_{L^2}^2 \\
   \geq&  \frac 12 \( 1-(\frac{\|X_0\|_{L^2}}{\|Q\|_{L^2}})^{\frac{4}{d}}\) \|\na (e^{i\a \phi}X(t))\|_{L^2}^2
   \geq 0.
\end{align*}
On the other hand, a straightforward expansion shows that
\begin{align*}
   H(e^{i\a \phi} X(t))
   = \frac{1}{2}\|\na \phi X(t)\|_{L^2}^2 \a^2 -  \a ({\rm Im} \int X(t) \na \overline{X(t)} \cdot \na \phi
   dx) + H(X(t)).
\end{align*}
Thus,  we arrive at
\begin{align}
   \frac{1}{2}\|\na \phi X(t)\|_{L^2}^2 \a^2 -  \a ({\rm Im} \int X(t) \na \overline{X(t)} \cdot \na \phi
   dx) + H(X(t)) \geq 0,
\end{align}
which yields \eqref{esti-integ}
due to the arbitrariness of $\a$. \hfill $\square$

\begin{lemma} [Control of Hamiltonian]  \label{Lem-bdd-H}
Let $X$ be the solution to \eqref{SNLS} on $[0,\tau^*)$.
Suppose that $\|X_0\|_{L^2} \leq \|Q\|_{L^2}$.
Then, for any $T>0$,
\begin{align} \label{bdd-H}
  \bbe \sup\limits_{t\in[0,\tau^*\wedge T)}   H(X(t))  \leq C(T) <\9.
\end{align}
In particular,
\begin{align} \label{bdd-H-a.s.}
   \sup\limits_{t\in[0,\tau^*\wedge T)}  H(X(t))  <\9,\ \ \bbp-a.s..
\end{align}
\end{lemma}

We remark that
unlike in the NLS case,
the Hamiltonian fails to be conserved in the stochastic case.
Lemma \ref{Lem-bdd-H} shows that
the finiteness of the Hamiltonian can be derived in every bounded time regime,
which is sufficient for the blow-up analysis.

{\it \bf Proof of Lemma \ref{Lem-bdd-H}.}
As in \cite[Theorem 3.1]{BRZ16} we have the evolution formula of the Hamiltonian
\begin{align} \label{H-evo}
   H(X(t)) &= H(X_0) + \sum\limits_{l=1}^N \int_0^t \frac{1}{2} \|\na   \phi_l X(s)\|_{L^2}^2 ds
              -  \sum\limits_{l=1}^N  \int_0^t {\rm Im} \<\na \phi_l X(s), \na X(s)\>_2 dB_l(s)  \nonumber \\
          & =: H(X_0) + H_1(t) + H_2(t).
\end{align}

Note that, by the conservation law of mass,
\begin{align*}
  \bbe \sup\limits_{s\in[0,\tau^*\wedge t)} H_1(s)
 \leq \frac{1}{2} (\sum\limits_{l=1}^N\|\na   \phi_l\|_{L^{\9}}^2) \|X_0\|_{L^2}^2 t.
\end{align*}
Moreover, by the Burkholder-Davis-Gundy inequality and
the Banica type inequality \eqref{esti-integ},
\begin{align*}
   \bbe \sup\limits_{s\in[0,\tau^*\wedge t)} |H_2(s)|
  \leq&C \bbe \(\int_0^{\tau^*\wedge t} \sum\limits_{l=1}^N |{\rm Im}\<\na\phi_l X(s), \na
  X(s)\>_2|^2ds\)^{\frac{1}{2}}\\
  \leq& \sqrt{2}C \bbe \(\int_0^{\tau^*\wedge t} \sum\limits_{l=1}^N H(X(s)) \|X(s)\na \phi_l\|_{L^2}^2
  ds\)^{\frac{1}{2}}\\
  \leq& \sqrt{2}C \|X_0\|_{L^2} \(\sum\limits_{l=1}^N\|\na\phi_l\|^2_{L^{\9}}\)^\frac 12
               \bbe\(\int_0^{\tau^*\wedge t} H(X(s))ds\)^{\frac{1}{2}}.
\end{align*}
Note that
$H(X(s))\geq 0$ for any $s<\tau^*$,
due to  Lemma \ref{Lem-GN} and the fact that $\|X(s)\|_{L^2} = \|X_0\|_{L^2} \leq \|Q\|_{L^2}$.

Thus, plugging the estimates above into \eqref{H-evo} we get that
\begin{align*}
  \bbe \sup\limits_{s\in[0, \tau^* \wedge t)} H(X(s))
  &\leq CT +C \bbe\(\int_0^t \sup\limits_{r\in [0,\tau^*\wedge s)} H(X(r))ds\)^{\frac{1}{2}}\\
  &\leq C(1+T) + C \int_0^t \bbe\sup\limits_{r\in[0,\tau^*\wedge s)} H(X(r))ds, \ \ t\in [0,T],
\end{align*}
which along with the Gronwall inequality implies that
\begin{align}
  \bbe\sup\limits_{t\in[0,\tau^*\wedge T)} H(X(t)) \leq C(1+T)e^{CT}<\9,
\end{align}
thereby yielding \eqref{bdd-H} and finishing the proof. \hfill $\square$

We are now ready to prove Theorem \ref{Thm-L2-Con}.

{\it \bf Proof of Theorem \ref{Thm-L2-Con}.}
Set
$\lbb(t) :=\frac{\|\na Q\|_{L^2}}{\|\na X(t)\|_{L^2}}$
and rescale
$X_{\lbb}(t,x) :=\lbb(t)^{\frac{d}{2}} X(t,\lbb(t)x)$,
$t\in [0,\tau^*)$, $x\in \bbr^d$.
Then, for any $t\in [0,\tau^*)$,
\begin{align} \label{scale-XQ}
\|X_{\lbb}(t)\|_{L^2}^2 = \|X(t)\|_{L^2}^2=\|Q\|_{L^2}^2, \ \
\|\na X_{\lbb}(t)\|_{L^2}^2= \lbb(t)^2 \|\na X(t)\|_{L^2}^2 = \|\na Q\|_{L^2}^2.
\end{align}

In view of \eqref{Hamil-def} and \eqref{H-evo},
we deduce that for $t\in [0,\tau^*)$,
\begin{align} \label{sca-H}
   H(X_{\lbb}(t))
   =& \lbb^2(t) H(X(t))  \nonumber \\
   =& \|\na Q\|_{L^2}^2  \(\frac{H(X_0)}{\|\na X(t)\|_{L^2}^2}
       +  \frac{1}{2} \sum\limits_{l=1}^N  \frac{\int_0^t\|\na \phi_l X(s)\|_{L^2}^2ds}{\|\na X(t)\|_{L^2}^2}
       - \frac{M(t)}{\|\na X(t)\|_{L^2}^2} \),
\end{align}
where
\begin{align*}
M(t) := \sum\limits_{l=1}^N \int_0^t {\rm Im} \<\na \phi_l X(s), \na X(s)\>_2dB_l(s).
\end{align*}

Note that since $\int_0^t \|\na \phi_l X(t)\|_{L^2}^2 ds \leq \|\na \phi_l\|^2_{L^\9} \|X_0\|^2_{L^2} t$
and $\|\na X(t)\|_{L^2} \to \9$ as $t\to \tau^*$,
\begin{align} \label{H12-0}
   \|\na Q\|_{L^2}^2  \(\frac{H(X_0)}{\|\na X(t)\|_{L^2}^2}
       + \frac{1}{2}\sum\limits_{l=1}^N  \frac{\int_0^t\|\na \phi_l X(s)\|_{L^2}^2ds}{\|\na X(t)\|_{L^2}^2} \)
   \to 0, \ \ as\ t\to \tau^*.
\end{align}

Regarding the  term $M(t)$,
using the time-change for martingales (see e.g. \cite[Theorem 16.4]{K97})
we deduce that there exists a Brownian motion $\wt B$
such that
\begin{align} \label{BM-M}
   M(t)= \wt B (\<M(t)\>), \ \ t\in (0,\tau^*),\  a.s.,
\end{align}
where $\<M\>$ denotes the quadratic variation process of $M$,
i.e.,
\begin{align} \label{M-quad}
    \<M\>_t = \sum\limits_{l=1}^N \int_0^t ( {\rm Im} \<\na \phi_l X(s), \na X(s)\>_2)^2 ds,\ \ t\in [0,\tau^*).
\end{align}
Moreover, by the law of iterated logarithm for Brownian motions (see \cite[Theorem 11.18]{KS91}),
\begin{align} \label{BM-lnln}
   \limsup\limits_{t\to \tau^*} \frac{|M(t)|}{\sqrt{2\<M(t)\>\ln\ln
   \<M(t)\>}} =1, \  \  a.s..
\end{align}

Hence,
we can take a measurable set $\Omega_0\subseteq \Omega$ with full probability
such that both \eqref{BM-M} and \eqref{BM-lnln} hold on $\Omega_0$.
Below we consider $\omega \in \{\tau^*<\9\} \cap \Omega_0$.
For simplicity,
the dependence on $\omega$ is omitted in the notation.

Let $t_n:= \inf\{t\in (0,\tau^*):  \|\na X(t)\|_{L^2} \geq n\}$,
$n\geq N_0:= 2 \|X_0\|_{H^1}$.
Then, using the continuity of $X$ in $H^1$ and
the explosion of $\|\na X(t)\|_{L^2}$ at $\tau^*$
we have
$t_n\to \tau^*$ and $\sup_{t\in[0,t_n]} \|\na X(t)\|_{L^2} = \|\na
X(t_n)\|_{L^2}$.

Then,
in the case where $\sup_{0\leq t<\tau^*}\<M(t)\><\9$,
using the continuity of Brownian motion
we have
$\sup_{0\leq t<\tau^*} |\wt B (\<M(t)\>)|<\9$,
and thus
\begin{align} \label{M.1}
   \frac{|M(t)|}{\|\na X(t)\|_{L^2}^2}
   \leq \frac{\sup\limits_{0\leq t< \tau^*} |\wt B (\<M(t)\>)|}{\|\na
   X(t)\|_{L^2}^2} \to 0,\ \ as~t\to \tau^*.
\end{align}

Moreover,
in the case where $\<M(t)\> \to \9$  as $t\to \tau^*$,
by \eqref{M-quad} and H\"older's inequality,
\begin{align*}
   \<M(t_n)\>
   \leq& \sum\limits_{l=1}^N \|\na \phi_l\|_{L^\9}^2 \|X_0\|_{L^2}^2 \int_0^{t_n} \|\na X(s)\|_{L^2}^2 ds \nonumber \\
   \leq& \sum\limits_{l=1}^N t_n  \|\na \phi_l\|_{L^\9}^2 \|X_0\|_{L^2}^2 \|\na X(t_n)\|_{L^2}^2 .
\end{align*}
which along with \eqref{BM-lnln} yields that
\begin{align} \label{M.2}
   \limsup\limits_{n\to \9} \frac{|M(t_n)|}{\|\na X(t_n)\|_{L^2}^2}
   \leq& \limsup\limits_{n\to \9}  \frac{\sqrt{2\<M(t_n)\> \ln\ln \<M(t_n)\>}}{\|\na X(t_n)\|_{L^2}^2} \nonumber  \\
   \leq& C\limsup\limits_{n\to \9} \frac{\sqrt{\ln\ln \|\na X(t_n)\|_{L^2}}}{\|\na X(t_n)\|_{L^2}}
    = 0.
\end{align}

Thus,
we conclude from \eqref{M.1} and \eqref{M.2} that
\begin{align*}
   \frac{M(t_n)}{\|\na X(t_n)\|_{L^2}^2} \to 0,\ \ as~n\to \9,
\end{align*}
Taking into account \eqref{sca-H} and \eqref{H12-0},
we obtain that the Hamiltonian of the rescaled solution
tends to zero along the sequence $\{t_n\}$, namely,
\begin{align} \label{H-0}
 \lim\limits_{n\to \9} H(X_{\lbb_n}(t_n))=0,
\end{align}
where $\lbb_n:=\lbb(t_n) = \|\na Q\|_{L^2}/ \|\na X(t_n)\|_{L^2}$, $n\geq 1$.

Now, it follows from \eqref{scale-XQ} and \eqref{H-0} that
\begin{align*}
   \lim\limits_{n\to\9} \frac{d}{4+2d}
     \|X_{\lbb_n}(t_n)\|_{L^{2+\frac{4}{d}}}^{2+\frac{4}{d}}
   = \lim\limits_{n\to\9} ( \frac{1}{2}\|\na X_{\lbb_n}(t_n)\|_{L^2}^2 + H(X_{\lbb_n}(t_n)) )
   =\frac{1}{2} \|\na Q\|_{L^2}^2,
\end{align*}
which implies that
\begin{align}
   \lim\limits_{n\to\9} \|X_{\lbb_n}(t_n)\|_{L^{2+\frac{4}{d}}}
   =(\frac{2+d}{d})^{\frac{d}{4+2d}} \|\na Q\|_{L^2}^{\frac{d}{2+d}}.
\end{align}

Thus,  by virtue of Theorem \ref{Thm-Compact},
we obtain  a profile $V\in H^1$ and a sequence of
concentration points $\{y_n\} \subseteq \mathbb{R}^d$ such that
up to a subsequence (still denoted by $\{n\}$)
\begin{align}
   & \|Q\|_{L^2} \leq \|V\|_{L^2} , \label{v-q} \\
   X_{\lbb_n} (t_n,\cdot + y_n) & \rightharpoonup V\ \ weakly\ in\ H^1, \ as\ n\to \9. \label{weak-x-v}
\end{align}
In particular, for  any $A>0$ and $R>0$,
\begin{align} \label{esti-V}
   \int\limits_{|x|\leq A} |V(x)|^2 dx
   \leq&  \liminf\limits_{n\to\9} \int\limits_{|x|\leq A} \lbb_n^d
   |X(t_n, \lbb_n x + y_n)|^2 dx  \nonumber \\
   =&  \liminf\limits_{n\to\9} \int\limits_{|x-y_n|\leq \lbb_nA}
   |X(t_n, x)|^2 dx \nonumber \\
   \leq& \liminf\limits_{n\to\9} \int\limits_{|x-y_n|\leq R}
   |X(t_n, x)|^2 dx,
\end{align}
where  the last step is due to the fact that
$\lbb_n A \leq R$ for $n$ large enough.

Therefore,
letting $A \to \9$
and using \eqref{v-q}
we obtain \eqref{esti-L2-conc}
and finish the proof.   \hfill $\square$

\subsection{Universality of stochastic critical mass blow-up profiles} \label{Subsec-Self}

We first recall the variational characterization of the ground state.

\begin{lemma} [Variational characterization of ground state, \cite{R13}] \label{Lem-Ground}
Let $v\in H^1$ be such that
\begin{align*}
   \|v\|_{L^2} = \|Q\|_{L^2}, \ \ H(v) =0.
\end{align*}
Then, there exist parameters $\lbb_0\in \bbr^+$, $x_0\in \bbr^d$ and $\theta_0\in \bbr$
such that
\begin{align*}
   v(x) = \lbb_0^{\frac d 2} Q(\lbb_0 x + x_0) e^{i\theta_0},\ \ x\in \bbr^d.
\end{align*}
\end{lemma}
\medskip

Now we are prepared to prove Theorem \ref{Thm-Self}.

{\it \bf Proof of Theorem \ref{Thm-Self}.}
Let $\Omega_0$, $t_n$, $y_n$ and the profile $V$ be as in Theorem \ref{Thm-L2-Con}.
Below we consider $\omega \in \{\tau^*<\9\} \cap \Omega_0$.
For simplicity,
we suppress the $\omega$ dependence
in the notation.

We claim that
\begin{align}
    & \|V\|_{L^2} =\|Q\|_{L^2}, \label{V-Q}  \\
    & \|\na V\|_{L^2} = \|\na Q\|_{L^2},  \label{naV-naQ.1} \\
    & H(V) =0. \label{H-V-0}
\end{align}

For this purpose, using \eqref{v-q} and the conservation of mass we note that
\begin{align} \label{esti-XV.1}
   \|X_{\lbb_n}(\cdot + y_n)\|_{L^2} =\|X_0\|_{L^2} =\|Q\|_{L^2} \leq \|V\|_{L^2},
\end{align}
where $\lbb_n := \lbb(t_n)$.
Moreover,  \eqref{esti-V} yields that,
\begin{align} \label{esti-XV.2}
   \|V\|_{L^2}^2 \leq  \liminf\limits_{n\to \9}
   \int\limits_{|x-y_n|\leq  R} |X(t_n, x)|^2 dx
   \leq \|X(t_n)\|_{L^2}^2 = \|Q\|_{L^2}^2.
\end{align}
Thus,  \eqref{V-Q} follows immediately from \eqref{esti-XV.1} and \eqref{esti-XV.2}.

Regarding \eqref{naV-naQ.1} and \eqref{H-V-0},
in view of \eqref{weak-x-v} and \eqref{esti-XV.1},
we first have
\begin{align} \label{Xlbb-V-L2}
   X_{\lbb_n}(t_n,\cdot + y_n) \to V\ \ in\ L^2, \ as\ n\to \9,
\end{align}
which along with the uniform boundedness of $\|\na X_{\lbb_n}(t_n)\|_{L^2}$,
implied by \eqref{scale-XQ}, and interpolation
yields the strong convergence  in $L^{2+4/d}$, i.e.,
\begin{align} \label{strong-X-V-Lp}
   X_{\lbb_n}(t_n,\cdot + y_n) \to V\ \ in\ L^{2+\frac{4}{d}},\ as\ n\to \9.
\end{align}

In order to prove the strong convergence in $H^1$,
i.e.,
\begin{align} \label{strong-X-V-H1}
     X_{\lbb_n}(t_n, \cdot + x_n) \to V\ \ in\ H^1,\ as\ n\to \9.
\end{align}
we note that, on one hand, by \eqref{scale-XQ}, \eqref{H-0} and \eqref{strong-X-V-Lp},
\begin{align*}
   \|\na Q\|_{L^2}^2
   =& \lim\limits_{n\to \9}\|\na X_{\lbb_n}(\cdot + y_n) \|_{L^2}^2 \\
   =& \lim\limits_{n\to \9}
       2  H(X_{\lbb_n}(\cdot + y_n))
       + \frac{d}{d+2} \|X_{\lbb_n}(\cdot + y_n) \|_{L^{2+\frac 4d}}^{2+\frac 4d} \\
   =& \frac{d}{d+2}\|V\|_{L^{2+\frac{4}{d}}}^{2+\frac 4 d},
\end{align*}
which along with  Lemma \ref{Lem-GN}
and \eqref{V-Q} yields
\begin{align*}
   \|\na Q\|_{L^2}^2 \leq \|\na V\|_{L^2}^2.
\end{align*}
On the other hand,
\eqref{scale-XQ} and \eqref{weak-x-v} yield that
\begin{align}  \label{naV-naQ}
   \|\na V\|_{L^2}^2 \leq \liminf\limits_{n\to \9} \|\na X_{\lbb_n}(t_n, \cdot + y_n)\|_{L^2}^2 = \|\na Q\|_{L^2}^2.
\end{align}
Thus, we obtain
$$\|\na V\|_{L^2}^2 = \|\na Q\|_{L^2}^2 =
\|\na X_{\lbb_n}(\cdot + y_n) \|_{L^2}^2, $$
which yields \eqref{naV-naQ.1},
and, via \eqref{weak-x-v} and \eqref{Xlbb-V-L2},
also gives the strong convergence \eqref{strong-X-V-H1} in $H^1$.

Therefore, combining \eqref{strong-X-V-Lp} and \eqref{strong-X-V-H1},
we obtain \eqref{H-V-0}, as claimed.

Now, by virtue of \eqref{V-Q}, \eqref{naV-naQ.1}, \eqref{H-V-0} and the
variational characterization  of the ground state $Q$ in Lemma \ref{Lem-Ground},
we obtain that there exist $x_0\in \bbr^d$ and $\theta_0\in [0,2\pi)$ such
that
\begin{align}
   V(x) =e^{i\theta_0} Q(x+x_0),
\end{align}
which along with \eqref{strong-X-V-H1} implies \eqref{Strong-X-Q}.
The proof is complete.
\hfill $\square$
\medskip

{\bf Proof of Corollary \ref{Cor-L2-conv}.}
In order to prove \eqref{X2-delta},
it suffices to prove that
for any $ \vf \in C_b^\9$,
\begin{align} \label{X2-vf-0}
    \int \vf(x) |X(t_n,x)|^2 dx - \|Q\|_{L^2}^2 \vf(y_n) \to 0,\ \ as\ n\to \9.
\end{align}

For this purpose,
by a change of variables
we obtain
\begin{align*}
    \int \vf(x) |X(t_n,x)|^2 dx
    =& \int \vf(\lbb_n x +y_n) |\lbb_n^{\frac d2} X(t_n, \lbb_n x +y_n)|^2 dx.
\end{align*}
This yields
\begin{align}
   & \bigg|\int \vf(x) |X(t_n,x)|^2 dx - \vf(y_n) \|Q\|_{L^2}^2 \bigg| \notag\\
   \leq& \int |\vf(\lbb_nx + y_n)| | |\lbb_n^\frac d2 X(t_n, \lbb_n x + y_n)|^2 - |Q(x)|^2| dx
         + \int| \vf(\lbb_nx+y_n) - \vf(y_n)| |Q(x)|^2 dx \notag \\
   \leq& 2\|Q\|_{L^2} \|\vf\|_{L^\infty} \|\lbb_n^\frac d2 X(t_n, \lbb_n x + y_n)e^{i\theta(t_n)} - Q(x)\|_{L^2}
         + \lbb_n \|\vf'\|_{L^\infty} \int |x| |Q(x)|^2 dx.
\end{align}
Thus, in view of \eqref{Strong-X-Q},
the right-hand side above converges to zero.
The proof is complete.
\hfill $\square$

\subsection{Vanishing of virial} \label{Subsec-Explosion}

We take the full probability set $\Omega_0$ such that
Theorem \ref{Thm-Self}, Corollary \ref{Cor-L2-conv} and Lemma \ref{Lem-bdd-H} hold on $\Omega_0$.
For every $\omega \in \Omega_0$,
let $\{t_n(\omega)\}$ and $\{y_n(\omega)\}$ be as in Theorem \ref{Thm-Self}.
Up to a subsequence (still denoted by $n$),
we may assume that $y_n(\omega)\to y^*(\omega)$,
$y^*(\omega)$ may be the Alexandrov point of $\bbr^d$.
In the case where $|y^*(\omega)|<\9$,
we set
\begin{align} \label{Vm-def}
   V_m(t)(\omega):=\int \theta_m(x-y^*(\omega)) |X(t)(\omega)|^2 dx,
\end{align}
where $\theta_m(x):=m^2 \theta(\frac{x}{m})$,
$\theta\in C_0^{\9}$ satisfies $\theta(x)=|x|^2$ for $x\in B(0,1)$,
$|\na \theta(x)|^2\leq C \theta (x)$ for $x\in \bbr^d$, $C>0$.

In order to analyze the virial functional \eqref{Vm-def},
we shall derive an evolution formula on a full probability set
in Lemma \ref{Lem-Vm-Ito} below.
The delicate fact is that,
it is unclear whether $y^*$ is measurable with respect to $\omega$,
so the application of It\^o's formula as in \cite{BRZ14.3} is not possible directly.
Below we show that
the evolution formula can be derived for $ V_m(t)$, 
but with  $y^*(\omega)$ replaced by any fixed $y\in \mathbb{R}^d$.
The point is that this formula holds for any $y\in \mathbb{R}^d$
on a universal full probability set $\Omega_1$,
but the null set $\Omega\setminus \Omega_1$ is independent of $y$.
Thus, one can eventually replace $y$ by $y^*(\omega)$ in the evolution formula.

The proof is based on a mollification procedure.
Let $\vf \in C_c^\infty$ be a real-valued non-negative function with unit integral,
$\vf_\ve(x) := \ve^{-d}\vf(\frac x \ve)$, $x\in \bbr^d$.
Set $h^\ve : = h\ast \vf_\ve$ for any locally integrable function $h$ mollified by $\vf_\ve$,
$f^\ve(X):= (|X|^\frac{4}{d}X)^\ve$
and $f(X):= |X|^\frac{4}{d}X$.
We have that $\bbp$-a.s. for any $t\in (0,\tau^*)$
and any $q\in (1,\infty)$,
as $\ve \to 0$,
\begin{align}  \label{molif-limit}
   X_0^\ve \to X_0\  in\ L^2, \ \
   X^\ve \to X\ \in\ L^q(0,t; H^1)\cap L^{2+\frac 4d}(0,t;L^{2+\frac 4d}),\ \
   f^\ve(X) \to f(X)\ in\ L^{\frac{2d+4}{d+4}}(0,t;L^{\frac{2d+4}{d+4}}).
\end{align}

\begin{lemma} \label{Lem-Mart-0}
Consider the situation as in Theorem \ref{Thm-Var-0}.
Then, for any sequence $\ve_n\to 0$ when $n\to \infty$,
there exists a further subsequence $\{\ve_{n_j}\}$
and a full probability set $\wt \Omega_0$,
such that for any $\omega\in \wt \Omega_0$ and any $y\in \bbr^d \cap \mathbb{Q}^d$,
\begin{align} \label{betaj-0}
  \sum\limits_{l=1}^N \bigg| \int_0^t \int  \theta_m(x-y) \ol{X}^{\ve_{n_j}} (s,x)((X\phi_l)^{\ve_{n_j}} - X^{\ve_{n_j}} \phi_l)(s,x) dx d B_l(s) \bigg|
  \to 0,\ \ as\ \ve_{n_j} \to 0,\
  t\in (0,\tau^*).
\end{align}
\end{lemma}

{\bf Proof.}
Fix $1\leq l\leq N$.
For any $y\in \mathbb{R}^d\cap \mathbb{Q}^d$ and any $M\geq 1$,
we use the Burkholder-Davis-Gundy inequality
and the conservation of mass to derive
\begin{align*}
   & \bbe \sup_{s\in [0,\tau^*\wedge M)} \bigg|\int_0^{s} \int\theta_m(x-y) \ol{X}^{\ve_n}(s,x) ((X\phi_l)^{\ve_n} - X^{\ve_n} \phi_l)(s,x) dx d B_l(s) \bigg|^2  \\
  =&  \bbe  \( \int_0^{\tau^*\wedge M}  \bigg|  \int \theta_m(x-y)\ol{X}^{\ve_n}(s,x) ((X\phi_l)^{\ve_n} - X^{\ve_n} \phi_l)(s,x) dx \bigg|^2  d s \)^{\frac 12} \\
 \lesssim& m^2  \|X_0\|_{L^2}  \bbe \(\int_0^{\tau^*\wedge M}\|(X\phi_l)^{\ve_n}(s) - X^{\ve_n}(s) \phi_l\|_{L^2}^2 ds \)^\frac 12.
\end{align*}
Since $\bbp$-a.s. for any $s\in (0,\tau^*\wedge M)$,
\begin{align*}
   \|(X\phi_l)^{\ve_n}(s) - X^{\ve_n}(s) \phi_l\|_{L^2}^2 \to 0,\ \ as\ \ve_n \to 0,
\end{align*}
and
\begin{align*}
    \|(X\phi_l)^{\ve_n}(s) - X^{\ve_n}(s) \phi_l\|_{L^2}^2
    \lesssim  \|\phi_l\|_{L^\infty}^2  \|X(s)\|_{L^2}^2
    \lesssim  \|\phi_l\|_{L^\infty}^2  \|X_0\|_{L^2}^2,
\end{align*}
the dominated convergence theorem gives
\begin{align*}
 \bbe \sup_{s\in [0,\tau^*\wedge M)}
 \bigg|\int_0^{s} \int \theta_m(x-y) \ol{X}^{\ve_n}(s,x) ((X\phi_l)^{\ve_n} - X^{\ve_n} \phi_l)(s,x) dx d B_l(s) \bigg|^2
 \to 0,\ \ as\ \ve_n\to 0,
\end{align*}
which implies that
\begin{align*}
  \sup_{s\in [0,\tau^*\wedge M)} \bigg|\int_0^{s} \int \theta_m(x-y)\ol{X}^{\ve_n}(s,x) ((X\phi_l)^{\ve_n} - X^{\ve_n} \phi_l)(s,x) dx d B_l(s) \bigg|^2
 \to 0\ \ in\ probability,\ as\ \ve_n\to 0.
\end{align*}
It follows that
for any $y\in \bbr^d \cap \mathbb{Q}^d$
and any $M\geq 1$, there exists a subsequence $\{\ve_{n_j}\}$
such that $\bbp$-a.s.,
\begin{align} \label{betaj-tauM-0}
  \sup_{s\in [0,\tau^*\wedge M)} \bigg|\int_0^{s} \int \theta_m(x-y)\ol{X}^{\ve_{n_j}}(s,x) ((X\phi_l)^{\ve_{n_j}} - X^{\ve_{n_j}} \phi_l)(s,x) dx d B_l(s) \bigg|^2
 \to 0,\ \  as\ \ve_{n_j}\to 0.
\end{align}
Using a diagonal argument
we can select a further subsequence, still denoted by $\{\ve_{n_j}\}$,
and a full probability set $\wt \Omega_0$,
such that for any $y\in \bbr^d \cap \mathbb{Q}^d$
and any $M\geq 1$,
the convergence \eqref{betaj-tauM-0} holds on $\wt \Omega_0$.
Taking into account that 
\begin{align*}
   \{t<\tau^*\} = \bigcup_{M\geq 1} \{t<\tau^*\wedge M\}, 
\end{align*}
we thus obtain \eqref{betaj-0} and finish the proof.
\hfill $\square$

The key evolution formula of the virial $V_m$ defined in \eqref{Vm-def}
is contained in Lemma \ref{Lem-Vm-Ito} below.

\begin{lemma} [Evolution of virial] \label{Lem-Vm-Ito}
Consider the situation as in Theorem \ref{Thm-Var-0}.
Suppose that $|y^*|<\infty$.
Then, there exists a measurable set $\Omega_1$ with full probability such that
for every $\omega\in \Omega_1$,
\begin{align} \label{Vm-Ito}
   V_m(t)(\omega)
   = V_m(0)
     - 2 {\rm Im} \int_0^t \<\na \theta_m (\cdot-y^*(\omega)) X(s)(\omega), \na X(s)(\omega)\> ds,\ \
     0\leq t<\tau^*(\omega).
\end{align}
\end{lemma}

{\bf Proof.}
Taking convolution of both sides of \eqref{SNLS} with the mollifier $\vf_\ve$
we derive from \eqref{SNLS} that for every $x\in \mathbb{R}^d$,
\begin{align*}
   X^\ve(t,x) = \int_0^t (i\Delta X^\ve + if^\ve(X)  - (\mu X)^\ve)(s,x) ds
                + \sum\limits_{l=1}^N \int_0^t i X \phi_l^\ve(s,x) d B_l(s), \ \
   t\in (0,\tau^*).
\end{align*}
This along with the product rule yields
\begin{align} \label{X2-Ito}
   |X^\ve(t,x)|^2
  & = |X^\ve(0,x)|^2
     + 2 {\rm Re} \int_0^t [\ol{X}^\ve (i\Delta X^\ve + if^\ve(X)  - (\mu X)^\ve)](s,x) ds  \notag  \\
    & + \sum\limits_{l=1}^N \int_0^t |(X\phi_l)^\ve(s,x)|^2 ds
     - 2 {\rm Im} \sum\limits_{l=1}^N \int_0^t  \ol{X}^\ve(s,x) (X \phi_l)^\ve(s,x) dB_l(s),\ t\in(0,\tau^*), \mathbb{P}-a.s..
\end{align}
Since both sides are continuous in $x$,
we can find a universal null set $\wt N$ such that
for every $\omega\in \Omega\setminus \wt N$,
\eqref{X2-Ito} holds for any $x\in \bbr^d$.

Thus,
we can take a full probability set $\Omega_1$ such that
Theorem \ref{Thm-Self}, Lemmas \ref{Lem-bdd-H} and \ref{Lem-Mart-0}
and \eqref{X2-Ito} hold on $\Omega_1$.
Below we consider $\omega\in \Omega_1$ and omit the argument $\omega$ in the notation for simplicity.

For any $y\in \bbr^d\cap \mathbb{Q}^d$, set
\begin{align*}
   V_{m,y}^\ve(t) : = \int \theta_m(x-y) |X^\ve(t,x)|^2 dx,\ \ t\in (0,\tau^*).
\end{align*}
Using \eqref{X2-Ito} we have the following identity on $\Omega_1$:
\begin{align} \label{Vm-I1-I4}
   V^\ve_{m,y}(t)
   =& \int \theta_m(x-y) |X^\ve(0,x)|^2  dx
       + \sum\limits_{l=1}^N \int  \theta_m(x-y)\int_0^t |(X\phi_l)^\ve(s,x)|^2 ds dx \notag \\
    & + 2 {\rm Re} \int  \theta_m(x-y) \int_0^t [\ol{X}^\ve (i\Delta X^\ve + if(X)^\ve  - (\mu X)^\ve)](s,x) ds dx \notag \\
    & - 2 \sum\limits_{l=1}^N {\rm Im}  \int  \theta_m(x-y) \int_0^t \ol{X}^\ve (s,x) (X \phi_l)^\ve(s,x) d B_l(s) dx. \notag  \\
   =&: \int \theta_m(x-y) |X^\ve(0,x)|^2  dx
       + I_{1,y}^\ve(t) + I_{2,y}^\ve(t) + I_{3,y}^\ve(t),\ \ t\in (0,\tau^*).
\end{align}

Since $\|\theta_m(\cdot-y)\|_{L^\infty} \lesssim m^2$,
we infer from \eqref{molif-limit} that
\begin{align}  \label{I1-limit}
   \int \theta_m(x-y) |X^\ve(0,x)|^2  dx \to \int \theta_m(x-y) |X(0,x)|^2  dx = V_m(0), \ \ as\ \ve \to 0.
\end{align}
Below we pass to the limit $\ve\to 0$ for each term $I^\ve_{j,y}(t)$, $1\leq j\leq 3$.

First, using \eqref{molif-limit} again we have
\begin{align}  \label{I2-limit}
    I_{1,y}^\ve(t) \to 2 \int  \theta_m(x-y) \int_0^t \mu(x) |X(s,x)|^2 ds dx
    =  2 \int_0^t  \int   \theta_m(x-y)\mu(x) |X(s,x)|^2 dx   ds,
\end{align}
where $t\in (0,\tau^*)$ and the last step is due to the Fubini theorem.

Regarding the second term $I^\ve_{2,y}(t)$,
using Fubini's theorem and the integration-by-parts formula we derive
\begin{align*}
   2 {\rm Re} \int  \theta_m(x-y) \int_0^t \ol{X}^\ve(s,x) i\Delta X^\ve(s,x)  ds dx
   = 2 {\rm Im} \int_0^t \int   \na X^\ve (s,x)  \cdot  \na \theta_m(x-y) \ol{X}^\ve(s,x)dx ds,
\end{align*}
which along with \eqref{molif-limit} and the fact
that $\|\na \theta_m(\cdot-y)\|_{L^\infty} \lesssim m$
yields that as $\ve \to 0$,
\begin{align} \label{I3.1-limit}
    2 {\rm Re} \int  \theta_m(x-y) \int_0^t \ol{X}^\ve(s,x) i\Delta X^\ve (s,x) ds dx
    \to 2 {\rm Im} \int_0^t \int  \na X(s,x) \cdot \na \theta_m(x-y) \ol{X}(s,x) dx ds.
\end{align}
Moreover,
it follows from \eqref{molif-limit} again that
\begin{align}  \label{I3.2-limit}
   2 {\rm Re} \int  \theta_m(x-y) \int_0^t  \ol{X}^\ve(s,x)  if(X)^\ve (s,x)  ds dx
   \to - 2 {\rm Im} \int  \theta_m(x-y) \int_0^t  \ol{X} (s,x) f(X) (s,x)  ds dx
   =0,
\end{align}
and
\begin{align} \label{I3.3-limit}
    2 {\rm Re} \int  \theta_m(x-y) \int_0^t  \ol{X}^\ve (s,x) (\mu X)^\ve(s,x)   ds dx
    \to 2 {\rm Re} \int_0^t   \int  \theta_m(x-y) \mu(x) |X(s,x)|^2  dx  ds.
\end{align}
Thus, it follows from \eqref{I3.1-limit}, \eqref{I3.2-limit} and \eqref{I3.3-limit} that
as $\ve \to 0$,
\begin{align} \label{I3-limit}
   I_{2,y}^\ve(t) \to& 2 {\rm Im} \int_0^t  \int \na X(s,x)  \cdot \na \theta_m(x-y) \ol{X}(s,x) dx ds \notag \\
                 & -2 {\rm Re}  \int_0^t   \int  \theta_m(x-y)\mu(x) |X(s,x)|^2 dx  ds, \ \ t\in (0,\tau^*).
\end{align}

Finally, for the last stochastic term $I_{3,y}^\ve(t)$.
Since
\begin{align*}
    {\rm Im} \int \theta_m(x-y) \ol{X}^\ve(s,x) X^\ve (s,x)\phi_l(x) dx =0,
\end{align*}
using the stochastic Fubini theorem
and Lemma \ref{Lem-Mart-0}
we obtain that there exists a sequence $\{\ve_n\}$
such that the following convergence holds on
a full probability set $\Omega_1$:
for any $y\in \mathbb{R}^d \cap \mathbb{Q}^d$ and any $t\in (0,\tau^*)$,
\begin{align}  \label{I4-limit}
   I_{3,y}^{\ve_n}(t)
   = -2 \sum\limits_{l=1}^N {\rm Im} \int_0^t \int  \theta_m(x-y) \ol{X}^{\ve_{n}} (s,x)((X\phi_l)^{\ve_{n}} - X^{\ve_{n}} \phi_l)(s,x) dx d B_l(s)
   \to 0,\ \ as\ \ve_n\to 0.
\end{align}

Therefore,
we conclude from \eqref{Vm-I1-I4}, \eqref{I1-limit}, \eqref{I2-limit},
\eqref{I3-limit} and \eqref{I4-limit} that
on the set $\Omega_1$,
as $\ve_n \to 0$, for any $y\in \mathbb{R}^d \cap \mathbb{Q}^d$ and any $t\in (0,\tau^*)$,
\begin{align*}
   V^{\ve_n}_{m,y} (t) \to \int \theta_m(x-y) |X(0,x)|^2 dx
    +  2 {\rm Im} \int_0^t \int \na X(s,x) \cdot \na \theta_m(x-y) \ol{X}(s,x)  dx ds.
\end{align*}
Taking into account
\begin{align*}
   V^{\ve_n}_{m,y} (t) \to \int \theta_m(x-y) |X(t,x)|^2 dx,\ \ t\in (0,\tau^*),
\end{align*}
we obtain that for every $\omega\in \Omega_1$,
it holds that for any $y\in \mathbb{R}^d \cap \mathbb{Q}^d$,
\begin{align} \label{Vmy-Ito}
   \int \theta_m(x-y) |X(t,x)(\omega)|^2 dx
   =&  \int \theta_m(x-y) |X(0,x)|^2 dx \notag \\
    & - 2 {\rm Im} \int_0^t \< \na \theta_m(\cdot-y) {X}(s)(\omega),  \na X(s)(\omega) \> ds,\ \ t\in (0,\tau^*(\omega)).
\end{align}

We see that both sides of \eqref{Vmy-Ito} are continuous in $y$.
Actually, for any $y_n \to y$, $n\to \infty$,
since $X(\omega)\in C([0,t];H^1)$, $t\in(0,\tau^*(\omega))$,
we have
\begin{align*}
   |{X}(\omega) \na \theta_m(x-y) \cdot \na X(\omega)|
   \leq C m |{X}(\omega) ||\na X(\omega)| \in L^1((0,t)\times \bbr^d),
\end{align*}
where the  constant $C$ is independent of $\{y_n\}$.
Taking into account $\na \theta_m(x-y_n) \to  \na \theta_m(x-y)$
and using the dominated convergence theorem
we get
\begin{align*}
    \int_0^t \< \na \theta_m(\cdot-y_n) {X}(s)(\omega),  \na X(s)(\omega) \> ds
    \to  \int_0^t \< \na \theta_m(\cdot-y) {X}(s)(\omega),  \na X(s)(\omega) \> ds,\ \ as\ n\to \infty.
\end{align*}
Similar arguments also apply to the other two terms in \eqref{Vmy-Ito}.

Thus, we obtain that, on the set $\Omega_1$,
\eqref{Vmy-Ito} holds for any $y\in \mathbb{R}^d$.
In particular,
the null set is independent of $y \in \mathbb{R}^d$.
Therefore, we can
replace $y$ by $y^*(\omega)$ in \eqref{Vmy-Ito},
where $\omega\in \Omega_1$,
and prove \eqref{Vm-Ito}.
\hfill $\square$

Next, we show the boundedness of the virial
of solutions,
based on the control of the Hamiltonian in Lemma \ref{Lem-bdd-H}.

\begin{lemma} [Control of virial] \label{Lem-Vm}
Consider the situation as in Theorem \ref{Thm-Var-0}.
Set
\begin{align*}
   V(t):=\int |x|^2 |X(t,x)|^2 dx,\ \ t\in [0,\tau^*).
\end{align*}
Then, there exists a full probability set $\Omega_1$
such that for every $\omega\in \Omega_1$,
\begin{align} \label{esti-Vt}
   \sup\limits_{t\in [0,\tau^*(\omega))} V(t)(\omega) <\9.
\end{align}
\end{lemma}

{\bf Proof.}
We take the full probability sets $\Omega_1$
as in the proof of Lemma \ref{Lem-Vm-Ito},
and so Theorem \ref{Thm-Self}, Corollary \ref{Cor-L2-conv} and Lemmas \ref{Lem-bdd-H} and \ref{Lem-Vm-Ito}
hold on $\Omega_1$.
Below we consider $\omega\in \Omega_1$
and omit the argument $\omega$.

Recall that $\{t_n\}$ and $\{y_n\}$ are as in Theorem \ref{Thm-Self}
and $y_n\to y^*$,
where $y^*$ may be the Alexandrov point of $\bbr^d$.
In the case where $|y^*|<\9$,
applying Lemma \ref{Lem-Vm-Ito} we have
\begin{align} \label{V-Ito}
   V_m(X(t))
   = V_m(X_0)
     - 2 {\rm Im} \int_0^t \<\na \theta_m (\cdot -y^*) X(s), \na X(s)\> ds,\ \
     0\leq t<\tau^*.
\end{align}

By  \eqref{esti-integ}, \eqref{bdd-H-a.s.}
and $|\na \theta_m|^2 \leq C\theta_m$,
\begin{align*}
   | \frac{d}{dt} V_m(t)|
    \leq   C (H(X(t)) \|\na \theta_m(x-y^*) X(t)\|_{L^2}^2)^\frac 12
    \leq   C(\tau^*)  V_m^\frac 12(t)  ,
\end{align*}
where $C(\tau^*)$ is independent of $m$.
This yields that
\begin{align} \label{esti-vm}
   |\sqrt{V_m(t)}-\sqrt{V_m(0)}| \leq  \frac{1}{2}C(\tau^*) t,\ \ t\in[0,\tau^*).
\end{align}
In particular, we have the following uniform bound for all $m,n\geq 1$,
\begin{align}
   \sqrt{V_m(0)} \leq \frac{1}{2}C(\tau^*) \tau^* + \sqrt{V_m(t_n)}.
\end{align}
Applying Corollary \ref{Cor-L2-conv}
we infer that
$V_m(t_n)\to 0$ as $n\to \9$,
and so the following uniform-in-$m$ boundedness holds
\begin{align} \label{esti-vm0}
   \sqrt{V_m(0)} \leq \frac 12 C (\tau^*) \tau^*<\9.
\end{align}
In particular,
since $\theta_m(x-y^*) \to |x-y^*|^2$ as $m\to \infty$,
Fatou's lemma yields
\begin{align*}
   \int|x-y^*|^2 |X_0|^2 dx \leq \frac 14 C(\tau^*)^2(\tau^*)^2 <\9.
\end{align*}

Thus, using  \eqref{esti-vm}, \eqref{esti-vm0},
the conservation of mass and Fatou's lemma again we obtain
\begin{align*}
    V(t)
   \leq& |y^*|^2 \|X_0\|_{L^2}^2
        + \int |x-y^*|^2 |X(t)|^2 dx  \nonumber \\
    \leq& |y^*|^2 \|X_0\|_{L^2}^2 + \liminf\limits_{m\to \9}V_m(t) \nonumber \\
    \leq& |y^*|^2 \|X_0\|_{L^2}^2 + \liminf\limits_{m\to \9} (\frac 12 C(\tau^*) \tau^* + \sqrt{V_m(0)})^2 \nonumber \\
    \leq& |y^*|^2 \|X_0\|_{L^2}^2 + (C(\tau^*)\tau^*)^2 <\9.
\end{align*}
thereby yielding  \eqref{esti-Vt} in the case where $|y^*|<\9$.

The proof for the case where
$y^*$ is the Alexandrov point of $\bbr^d$ is similar.
In this case,
we modify
$V_m(t)$ by $\wt V_m(t):= \int \theta_m(x) |X(t,x)|^2dx$.
Then,
estimate \eqref{esti-vm} still holds
for $\wt V_m(t)$,
and by Corollary \ref{Cor-L2-conv},
$\wt V_m(t_n)\to 0$ as $n\to \9$.
Thus,
using similar arguments as those below \eqref{esti-vm}
we obtain \eqref{esti-Vt} and finish the proof.  \hfill $\square$
\medskip

{\it \bf Proof of Theorem \ref{Thm-Var-0}.}
Below we perform the analysis on the full probability set $\Omega_1$
defined in the proof of Lemma \ref{Lem-Vm-Ito}.
Recall that $\{y_n\}$ are the concentration points as in Corollary
\ref{Cor-L2-conv}
and
$y_n\to y^*$.

We claim that
\begin{align} \label{x*}
   y^*(\omega)=\frac{\lim\limits_{n\to \9}\int x |X(t_n(\omega),x)(\omega)|^2 d x}{\|Q\|_{L^2}^2}.
\end{align}

To this end,
suppose that
$\{y_n\}$ converges to infinity,
then for $R$ large enough and for $n$ very large,
$\{x:|x-y_n|\leq 1\} \subseteq \{|x| \geq R\}$,
which, via the $L^2$ concentration in Corollary \ref{Cor-L2-conv},
yields that
\begin{align*}
   \|Q\|_{L^2}^2
   \leq \int_{|x-y_n|\leq 1} |X(t_n, x)|^2 dx
   \leq \int_{|x|>R} |X(t_n, x)|^2 dx.
\end{align*}
But on the other hand,
by the finiteness of the virial \eqref{esti-Vt},
\begin{align}
    \int_{|x|>R} |X(t_n, x)|^2 dx
    \leq \frac{V(t_n)}{R^2}
    \leq \frac{\sup_{t\in[0,\tau^*)}V(t)}{R^2}
    \leq \ve
\end{align}
for $R$ large enough, which leads to a contradiction.
It follows that $|y^*|<\infty$.

Hence, we may take $R$ large enough such that $\{y_n\}\subseteq B_R(0):=\{x:|x|\leq R\}$.
Let $\chi_R \in C_0^\infty$ be a cut-off function
such that $\chi_R(x) =1$ for $|x|\leq R$, and $\chi_R(x)=0$ for $|x|\geq R+1$.
In particular, $\chi_R(y_n) y_n = y_n$ for all $n\geq 1$.
Then, we estimate
\begin{align*}
     \bigg|\int |X(t_n)|^2 x dx  - \|Q\|_{L^2}^2 y_n\bigg|
   \leq& \bigg|\int |X(t_n)|^2 (1-\chi_R(x))xdx \bigg|
         + \bigg|\int |X(t_n)|^2 \chi_R(x) x dx - \chi_R(y_n) y_n \|Q\|_{L^2}^2 \bigg|  \notag \\
   \leq& \frac{1}{R} {\sup_{t\in [0,\tau^*)}V(t)}
           + \bigg|\int |X(t_n)|^2 \chi_R(x) x dx - \chi_R(y_n) y_n \|Q\|_{L^2}^2 \bigg|.
\end{align*}
In view of Corollary \ref{Cor-L2-conv},
the second term above tends to zero,
which yields that
\begin{align}
   \limsup\limits_{n\to \infty} \bigg|\int |X(t_n)|^2 x dx  - \|Q\|_{L^2}^2 x_n  \bigg|
   \leq  \frac{1}{R} {\sup_{t\in [0,\tau^*)}V(t)} .
\end{align}
Then, using \eqref{esti-Vt} and letting $R \to \infty$
we obtain
\begin{align} \label{con-Xn-Q}
   \int |X(t_n)|^2 x dx - \|Q\|_{L^2}^2 y_n \to 0,\ \ as\ n\to \9.
\end{align}

Moreover, since $\{X(t)\}$ is a continuous and adapted process in $H^1$,
arguing as in the proof of \cite[Lemma 4.2]{BRZ14.3}
we have that for any $n,m\geq 1$,
\begin{align*}
  \bigg| \int |X(t_m)|^2 x dx - \int |X(t_n)|^2 x dx \bigg|
   \leq&   2 \int_{t_n}^{t_m} \bigg|{\rm Im} \int X(s) \na \overline{X(s)} dx\bigg| ds \notag \\
   \leq& C \int_{t_n}^{t_m} \(H(X(s))\|X(s)\|_{L^2}^2\)^\frac 12 ds \notag \\
   \leq&  C(\tau^*)  \|X_0\|_{L^2} |t_m-t_n| \to 0,
\end{align*}
where the last two steps follow from Lemmas \ref{Lem-integ} and \ref{Lem-bdd-H}.
This yields that the limit $\lim_{n\to \9} \int |X(t_n)|^2 x dx$ exists.
Thus,
taking into account \eqref{con-Xn-Q} we obtain \eqref{x*}, as claimed.

Now, let $ {V}_m(t) :=\int \theta_m(x-y^*) |X(t,x)|^2 dx$.
Arguing as in the proof of \eqref{esti-vm} we obtain
that for any $t_n \not= t_l<\tau^*$,
\begin{align} \label{til-Vm}
   |\sqrt{{V}_m(t_{n})}-\sqrt{ {V}_m(t_{l})}|
   \leq  \frac{1}{2} C(\tau^*) |t_{n}-t_{l}|.
\end{align}
Since  by Corollary \ref{Cor-L2-conv},
$ {V}_m(t_{l}) \to 0$
as $l\to \9$,
we infer that
\begin{align}
   \sqrt{ {V}_m(t_{n})} \leq \frac{1}{2} C(\tau^*) (\tau^*-t_{n})
\end{align}
with $ C(\tau^*)$ independent of $m$ and $n$.
Then,   Fatou's lemma yields that as $n\to \9$,
\begin{align} \label{esti-V-T}
  \int|x-y^*|^2 |X(t_{n})|^2 dx
  \leq
  \liminf\limits_{m\to\9} {V}_m(X(t_{n}))
  \leq \frac{1}{4}C(\tau^*)^2 (\tau^*-t_{n})^2
  \to 0,
\end{align}
thereby yielding \eqref{conv-x-x*}.

Furthermore,
by the uncertainty principle,
\begin{align*}
   \int |X(t_n, x+y^*)|^2 dx
   \leq C \int |x|^2  |X(t_n, x+y^*)|^2 dx
          \int |\na X(t_n, x+y^*)|^2 dx,
\end{align*}
which along with the conservation of mass
and \eqref{esti-V-T} implies that
\begin{align*}
   \|Q\|_{L^2}^2
   = \int |X(t_n, x)|^2 dx
   \leq& C \int |x-y^*|^2 |X(t_n, x)|^2  dx
          \int |\na X(t_n, x)|^2 dx \\
   \leq& \frac{1}{4}CC^2(\tau^*) (\tau^*-t_n)^2 \|\na X(t_n)\|_{L^2}^2,
\end{align*}
thereby yielding \eqref{blowup-rate}.
Therefore, the proof of Theorem \ref{Thm-Var-0} is complete.
\hfill $\square$

\section*{Acknowledgements}

The authors thank for the financial support
by the Deutsche Forschungsgemeinschaft
(DFG, German Science Foundation) through SFB 1283/2
2021-317210226 at Bielefeld University.
D. Zhang  is also grateful for the support by NSFC (No. 12271352)
and Shanghai Rising-Star Program 21QA1404500.

\end{document}